\documentclass[12pt,leqno]{article}
\usepackage{amsfonts}
\usepackage{graphics}
\usepackage{amsmath}
\usepackage{rlepsf}

\newtheorem{Theorem}{Theorem}

\newtheorem{Corollary}[Theorem]{Corollary}

\newtheorem{Definition}[Theorem]{Definition}

\newtheorem{Lemma}[Theorem]{Lemma}
\newtheorem{Proposition}[Theorem]{Proposition}

\newtheorem{Fundamental Theorem}{Fundamental Theorem}

\newenvironment{Proof}[1][Proof]{\textbf{#1.} }{\ \rule{0.5em}{0.5em}}

\def \cD {\mathcal{D}}
\def \fin {\mathrm{fin}}
\def \K {\mathcal{K}}

\def \cP {\mathcal{P}}

\def \aa {\underline{\a}}
\def \ii {\underline{i}}
\def \jj {\underline{j}}
\def \ro{\rho \hskip-.06in\raise.09in \hbox{$\scriptstyle{o}$}}
\def \VgM{{V\hskip-.17in\raise.15in \hbox{$\scriptstyle{\g+1}$}}}
\def \Vg{V\hskip-.11in\raise.15in \hbox{$\scriptstyle{\g}$}}
\def \Vgm{V\hskip-.17in\raise.15in \hbox{$\scriptstyle{\g-1}$}}
\def \gh{g \hskip-.10in\raise.13in \hbox{$\scriptstyle{\frac{1}{2}}$}}

\def \C {\mathbb{C}}
\def \R{\mathbb{R}}
\def \Ch {\mathbb{C}[[h]]}
\def \N {\mathbb{N}}
\def \Nz {\mathbb{N}_0}
\def \hN {\frac{1}{2}\mathbb{N}_0}
\def \su {\mathfrak{su}(2)}
\def \sl {\mathfrak{sl}(2,\C)}
\def \ia {{i_\a}}
\def \ja {{j_\a}}
\def \ib {{i_\b}}
\def \jb {{j_\b}}
\def \ig {i_\g}

\def \sp {\mathrm{span}}

\def \W {\Omega}

\def \Ct {\mathcal{C}}

\def \slc {\mathfrak{sl}(2,\C)}
\def \Uslc {U(\mathfrak{sl}(2,\C))}

\def \Z {\mathbb{Z}}
\def \cZ {\mathcal{Z}}
\def \H {\mathcal{H}}
\newcommand\CGL [6] {\left ( \begin{matrix}#4 & #5\\#1 & #2\end{matrix}
\Big | \begin{matrix} #3\\#6 \end{matrix} \right )}

\newcommand\CGR [6] {\left ( \begin{matrix}#4\\ #1\end{matrix} \Big |
\begin{matrix} #2 &#3\\#5 &#6 \end{matrix} \right )}

\def \h {[[h]]}
\def \tn {\otimes}
\def \id {\mathrm{id}}
\def \op {\oplus}
\def \ctn {\hat{\otimes}}
\def \A {\mathcal{A}}
\def \su {U_q(\mathfrak{su}(2))}
\def \pol{\mathrm{Pol}(SU_q(2))}
\def \F{\mathcal{F}}
\def \Z{\mathbb{Z}}
\def \D {\Delta}

\newcommand \La [4] {\Lambda ^{#2 #3}_{#1 #4}}
\def \t {\tau}

\def \w {\omega}
\def \r {\rho}
\def \d {\delta}
\def \f {\phi}

\def \a {\alpha}
\def \l {\lambda}
\def \s {\sigma}
\def \b {\beta}
\def \g {\gamma}

\def \lg {\mathfrak{g}}

\def \ua{u\hskip-.08in\raise.09in \hbox{$\scriptstyle{\a}$}}
\def \uai{u\hskip-.08in\raise.09in \hbox{$\scriptstyle{\a}$}_{i_\a}}

\def \ug{u\hskip-.08in\raise.09in \hbox{$\scriptstyle{\g}$}}
\def \ugi{u\hskip-.08in\raise.09in \hbox{$\scriptstyle{\g}$}_{i_\g}}
\def \va{v\hskip-.08in\raise.09in \hbox{$\scriptstyle{\a}$}}
\def \vai{v\hskip-.08in\raise.09in \hbox{$\scriptstyle{\a}$}_{i_\a}}
\def \vaj{v\hskip-.08in\raise.09in \hbox{$\scriptstyle{\a}$}_{j_\a}}
\def \vaa{v\hskip-.08in\raise.09in \hbox{$\scriptstyle{\a_1}$}}
\def \van{v\hskip-.08in\raise.09in \hbox{$\scriptstyle{\a_n}$}}

\def \va{v\hskip-.08in\raise.09in \hbox{$\scriptstyle{\a}$}}
\def \vai{v\hskip-.08in\raise.09in \hbox{$\scriptstyle{\a}$}_{i_\a}}
\def \vaj{v\hskip-.08in\raise.09in \hbox{$\scriptstyle{\a}$}_{j_\a}}
\def \vaa{v\hskip-.08in\raise.09in \hbox{$\scriptstyle{\a_1}$}}
\def \van{v\hskip-.08in\raise.09in \hbox{$\scriptstyle{\a_n}$}}

\def \vb{v\hskip-.08in\raise.09in \hbox{$\scriptstyle{\b}$}}
\def \vbi{v\hskip-.08in\raise.09in \hbox{$\scriptstyle{\b}$}_{i_\b}}
\def \vbj{v\hskip-.08in\raise.09in \hbox{$\scriptstyle{\b}$}_{j_\b}}
\def \vbb{v\hskip-.08in\raise.09in \hbox{$\scriptstyle{\b_1}$}}
\def \vbn{v\hskip-.08in\raise.09in \hbox{$\scriptstyle{\b_n}$}}

\def \vg{v\hskip-.08in\raise.09in \hbox{$\scriptstyle{\g}$}}
\def \vgi{v\hskip-.08in\raise.09in \hbox{$\scriptstyle{\g}$}_{i_\g}}
\def \vgj{v\hskip-.08in\raise.09in \hbox{$\scriptstyle{\g}$}_{j_\g}}
\def \vgg{v\hskip-.08in\raise.09in \hbox{$\scriptstyle{\g_1}$}}
\def \vgn{v\hskip-.08in\raise.09in \hbox{$\scriptstyle{\g_n}$}}

\def \vd{v\hskip-.08in\raise.09in \hbox{$\scriptstyle{\d}$}}
\def \vdi{v\hskip-.08in\raise.09in \hbox{$\scriptstyle{\d}$}_{i_\d}}
\def \vdj{v\hskip-.08in\raise.09in \hbox{$\scriptstyle{\d}$}_{j_\d}}
\def \vdd{v\hskip-.08in\raise.09in \hbox{$\scriptstyle{\d_1}$}}
\def \vdn{v\hskip-.08in\raise.09in \hbox{$\scriptstyle{\d_n}$}}

\def \vo{v\hskip-.06in\raise.09in \hbox{$\scriptstyle{0}$}_0}
\def \voo{v\hskip-.06in\raise.09in \hbox{$\scriptstyle{0}$}^0}
\def \vaM{v\hskip-.12in\raise.09in \hbox{$\scriptstyle{\a+1}$}}
\def \vam{v\hskip-.12in\raise.09in \hbox{$\scriptstyle{\a-1}$}}

\def \ra{\rho \hskip-.08in\raise.09in \hbox{$\scriptstyle{\a}$}}
\def \Va{V\hskip-.1in\raise.15in \hbox{$\scriptstyle{\a}$}}

\def \Vz{V\hskip-.1in\raise.15in \hbox{$\scriptstyle{z}$}}

\def \Vo{V\hskip-.1in\raise.15in \hbox{$\scriptstyle{0}$}}

\def \Vu{V\hskip-.1in\raise.15in \hbox{$\scriptstyle{1}$}}

\def \Vpm{V\hskip-.15in\raise.15in \hbox{$\scriptstyle{p-1}$}}

\def \Xa {X\hskip-.1in\raise.15in \hbox{$\scriptstyle{\a}$}}
\def \Xai {{X\hskip-.1in\raise.15in
\hbox{$\scriptstyle{\a}$}}^{i_\a}_{j_\a}}
\def \SXai {q^{-i_\a+j_\a}(-1)^{2\a+i_\a-j_\a}{X\hskip-.1in\raise.15in
\hbox{$\scriptstyle{\a}$}}^{-j_\a}_{-i_\a}}

\def \sXai {X\hskip-.1in\raise.15in
\hbox{$\scriptstyle{\a}$}^{-j_\a}_{-i_\a}}

\def \sXbi {X\hskip-.1in\raise.15in
\hbox{$\scriptstyle{\b}$}^{-j_\b}_{-i_\b}}

\def \Xb {X\hskip-.1in\raise.15in \hbox{$\scriptstyle{\b}$}}

\def \Xak {X\hskip-.1in\raise.15in \hbox{$\scriptstyle{\a_k}$}}

\def \Xbi {{X\hskip-.1in\raise.15in
\hbox{$\scriptstyle{\b}$}}^{i_\b}_{j_\b}}
\def \SXbi {q^{-i_\b+j_\b}(-1)^{2\b+i_\b-j_\b}{X\hskip-.1in\raise.15in
\hbox{$\scriptstyle{\b}$}}^{-j_\b}_{-i_\b}}

\def \Xg {X\hskip-.1in\raise.15in \hbox{$\scriptstyle{\g}$}}
\def \Xgi {{X\hskip-.1in\raise.15in
\hbox{$\scriptstyle{\g}$}}^{i_\g}_{j_\g}}
\def \SXgi {q^{-i_\g+j_\g}(-1)^{2\g+i_\g-j_\g}{X\hskip-.1in\raise.15in
\hbox{$\scriptstyle{\g}$}}^{-j_\g}_{-i_\g}}

\def \Xgmi {{X\hskip-.1in\raise.15in
\hbox{$\scriptstyle{\g}$}}^{-j_\g}_{-i_\g}}

\def \ra{\rho \hskip-.06in\raise.09in \hbox{$\scriptstyle{\a}$}}

\def \rg{\rho \hskip-.06in\raise.09in \hbox{$\scriptstyle{\g}$}}

\def \rz{\rho \hskip-.06in\raise.09in \hbox{$\scriptstyle{z}$}}
\def \ro {\rho \hskip-.06in\raise.09in \hbox{$\scriptstyle{0}$}}

\def \rb{\rho \hskip-.06in\raise.09in \hbox{$\scriptstyle{\b}$}}

\def \ga{g \hskip-.08in\raise.09in \hbox{$\scriptstyle{\a}$}}
\def \gai{{g \hskip-.08in\raise.09in
\hbox{$\scriptstyle{\a}$}}^{i_\a}_{j_\a}}

\def \gaj{{g \hskip-.08in\raise.09in
\hbox{$\scriptstyle{\a}$}}^{j_\a}_{i_\a}}

\def \gak{{g \hskip-.08in\raise.09in \hbox{$\scriptstyle{\a_k}$}}}

\def \gbj{{g \hskip-.08in\raise.09in
\hbox{$\scriptstyle{\b}$}}^{j_\b}_{i_\b}}

\def \ggj{{g \hskip-.08in\raise.09in
\hbox{$\scriptstyle{\g}$}}^{j_\g}_{i_\g}}

\def \go{{g \hskip-.08in\raise.09in \hbox{$\scriptstyle{0}$}}^{0}_{0}}

\def \gb{g \hskip-.08in\raise.09in \hbox{$\scriptstyle{\b}$}}
\def \gbi{{g \hskip-.08in\raise.09in
\hbox{$\scriptstyle{\b}$}}^{i_\b}_{j_\b}}

\def \gg{g \hskip-.08in\raise.09in \hbox{$\scriptstyle{\g}$}}
\def \ggi{{g \hskip-.08in\raise.09in
\hbox{$\scriptstyle{\g}$}}^{i_\g}_{j_\g}}

\def \ctn {\hat{\tn}}

\begin{document}
\title{\  Knot Theory With The Lorentz Group}

\author{ Jo\~{a}o  Faria Martins\\ \footnotesize\it  {Departamento de Matem\'{a}tica, Instituto Superior T\'{e}cnico,}\\ {\footnotesize\it Av. Rovisco Pais, 1049-001 Lisboa, Portugal}\\ {\footnotesize\it jmartins@math.ist.utl.pt}}

\date{\empty}
\maketitle

\begin{abstract}
We analyse the perturbative expansion of the knot invariants defined from the unitary representations of the Quantum Lorentz Group in two different ways, namely using the Kontsevich Integral and weight systems, and the $R$-matrix in the Quantum Lorentz Group defined by Buffenoir and Roche. The two formulations are proved to be equivalent; and they both yield $\Ch$-valued knot invariants related with the Melvin-Morton expansion of the Coloured Jones Polynomial.                     
 \end{abstract}
\centerline{\textit{ 2000 Mathematics Subject Classification}: 57M27, 17B37, 20G42  } %Invariants of Knots and 3-Manifolds, Quantum groups (quantized enveloping algebras) and related deformations, Quantum groups (quantized function algebras) and their representations 

\section*{Introduction}

The main aim of this article is to show what a possible path to define knot invariants out of the infinite dimensional representations of the Lorentz Group is.

Let $\A$ be a Hopf algebra, its category of finite dimensional
representations is therefore a compact monoidal category. Let $q$ be a complex number not equal to $1$ or $-1$. Suppose
$\A=U_q(\lg)$ is the Drinfeld Jimbo algebra attached to the semisimple Lie
algebra $\lg$. Even though $\A$ is not a ribbon Hopf algebra, it
possesses a formal $R$-matrix and a formal ribbon element. These elements make
sense when applied to finite dimensional representations of $\A$, and thus its
category of finite dimensional representations is a ribbon category. This
means we have a knot invariant attached to any finite dimensional
representation of $\A$. This kind of knot invariants take values in $\C$.

A similar situation happens in the case of the Quantum Lorentz Group $\cD$ as defined by Woronovicz and Podle\'{s} in \cite{PW}. We shall use   especially the further developments in its theory by
Buffenoir and Roche, see  \cite{BR1} and \cite{BR2}. Despite the fact $\cD$ is not a Drinfeld Jimbo algebra, its structure of a quantum double, namely $\cD=\cD(\su,\mathrm{Pol}(SU_q(2))$ with $q \in (0,1)$, makes possible
the definition of a formal $R$-matrix on it. Also, it is  possible to define
a heuristic ribbon element. The category of finite dimensional
representations of $\cD$ can be proved to be a ribbon category, and thus we
can define knot invariants out of it.  In fact, as observed in \cite{BR2}, it is possible to prove that it is ribbon equivalent to the category of finite dimensional representations of $U_q(\mathfrak{su}(2) \tn_{R^{-1}} U_{q}(\mathfrak{su}(2))$. This last bialgebra equals $\su\tn \su$ as an algebra but has a coproduct twisted by $R^{-1}$, the inverse of the $R$-matrix of $\su$. This equivalence relates the knot invariants obtained  with  the Coloured Jones Polynomial in a nice way.

Such splitting of $\cD$ is not, however,  the most natural when considering unitary infinite dimensional representations of it. The general classification of the unitary representations of the Quantum Lorentz Group is due to Pusz, cf \cite{P}. In this case, as well as in the case of harmonic analysis, its definition as a quantum double is usually easier to deal with. A fact observed in
\cite{BR1} is that it is possible to describe the action of the
formal $R$-matrix of the Quantum Lorentz Group in a class of infinite dimensional representations of it. For this reason, it is  natural to ask whether there exists a knot
theory attached to the infinite representations of $\cD$. See also \cite{G}. We shall see the answer is affirmative at least in the perturbative level. Since we are working
with infinite dimensional representations the general formulation  of
Reshetikhin and Turaev for constructing knot invariants cannot be directly
applied. It is possible, though, given a knot diagram, or to be more
precise a  connected $(1,1)$-tangle diagram, to make a heuristic evaluation
of the Reshetikhin-Turaev functor on it. This  yields an infinite series for
any knot diagram. This method was also elucidated in  \cite{NR}. Unfortunately, at least for unitary infinite dimensional representations, these infinite series do not seem to converge at least for some simple knot diagrams. However, they  converge $h$-adicaly for $q=\exp(h/2)$, since the expansions of their terms as power series in $h$ starts increasing in degree. Therefore these evaluations do define $\Ch$-valued knot invariants. This article aims to define these invariants  from the Kontsevich Integral and weight systems. 

Let $\lg$ be a semisimple Lie algebra. The $h$-adic variant of Drinfeld Jimbo algebras, that is the algebras $U_h(\lg)$, is usually more practical to deal with if one wants to define
knot invariants out of the infinite representations $\lg$. Let us be given a knot $K$. The fact that $U_h(\lg)$ is a ribbon Hopf
algebra, and not merely a formal ribbon Hopf algebra,  makes it possible that a central element of it can be defined out of $K$, or to be more precise out of a $2$-dimensional diagram for it. See for example \cite{LM}. This central element is well defined and is a knot invariant. The centre of $U_h(\lg)$ is isomorphic, through a canonical isomorphism,  with the algebra of formal power series on the centre of the universal enveloping algebra $U(\lg)$ of $\lg$. This means that given a Lie algebra $\lg$ we have a knot invariant taking its values on the algebra of formal power series over the centre of $U(\lg)$. This invariant can be described out of the Kontsevich Integral. If we have an irreducible finite dimensional
representation $V$ of $\lg$ each of the terms of the formal power series
associated with the knot $K$ will then act in $V$ as a multiple of the
identity. Therefore we can transform a formal power series on the centre of
$U(\lg)$ into a formal power series over $\C$. If the representations are finite dimensional, these power series have a non zero radius of convergence and their
value at $h=2\log(q)$ is the value of the (rescaled) knot invariant associated with 
$U_q(\lg)$, as long as we use the representation of $U_q(\lg)$ that quantises
the representation $V$ of $\lg$ with which we are working. 

Notice that nothing says that the same framework cannot be applied to an infinite dimensional representations of $\lg$, as long as any central element of $U(\lg)$ acts in it as a multiple of the identity. Representations of this kind appear frequently in Lie algebra theory, and are commonly known as representations which admit a central character. Some examples are the irreducible cyclic highest weight representations of $\lg$, for $\lg$ semisimple, which, in the  $\slc$ context,  are simply constructed by perturbing the spin representation in such a way that we admit arbitrary complex spins. In this case this yields a knot invariant which is in some sense an analytic continuation of the Coloured Jones Polynomial.

Other examples of infinite dimensional representations that admit  a central character are the
representations of the Lie algebra $L$ of the Lorentz Group which correspond to
the representations of the Lorentz Group in the principal series. These are the
classical counterpart of the representations of the Quantum Lorentz
Group considered in \cite{BR2}. Therefore we would expect the knot power series invariants that come out of their use to  relate somehow with the knot invariants that come from the infinite dimensional representations of the Quantum Lorentz Group. A main result of this article is that the answer is yes.  Of course the way we construct central elements of the enveloping
algebra of $L$ must be specified. Notice that the Quantum Lorentz Group is not
the Drinfeld-Jimbo algebra associated with the Lorentz Algebra $L$. For the
reasons pointed out before, one solution is to  define the $h$-adic quantised
universal enveloping algebra of $L$ in a non standard  way as
$U_h(\mathfrak{su}(2)) \tn_{R^{-1}} U_{h}(\mathfrak{su}(2))$. Another
solution, which is equivalent, is to use the Kontsevich Universal Knot Invariant. Using it, we can associate  to a knot a series in the centre of the universal enveloping
algebra of $L$, as long as we specify an $L$ invariant, non degenerate, symmetric bilinear form in $L$. These series only depend on the  knot isotopy class. Such a bilinear form can be chosen so that the construction of central elements is coherent with the algebraic structure of
the Quantum Lorentz Group. The algebraic properties of this kind of  knot invariants will be a main topic of this article. 

We are mainly interested in the definition of numerical, rather than perturbative, knot invariants from infinite dimensional representations of the Quantum Lorentz Group. We expect our expansions to relate with them, if we can define any, as their perturbation series at the origin. These issues will be dealt with in a separate work, namely \cite{FM}, where the convergence properties of the power series obtained is analysed. A major result therein is that even though the power series can have a zero radius of convergence, they are, at least in some cases Borel-Gevrey summable. This indicates that some precise numerical knot invariants may be defined.

I finish referring to the main motivation of this work, namely its possible applications to Quantum Gravity. For an example of the use of the unitary representations of the Lorentz Group in the construction of spin foam models for  Quantum General Relativity we refer to \cite{BC}. See also \cite{NR} for its quantised counterpart.

\tableofcontents
\newpage

\section{Preliminaries}
\subsection{Chord Diagrams}

\begin{figure}
\begin{center}
\includegraphics{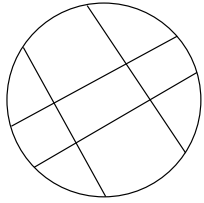}
\caption{\label{chord} A Chord Diagram with Four Chords.}
\end{center}
\end{figure}

We recall the definition  of the algebra of  chord diagrams, which is the
target space for the Kontsevich Universal Knot Invariant. For more details see
for example \cite {BN} or \cite{K}.
A chord diagram is a finite set  $w=\{c_1,...,c_n\}$ of cardinality two, non-intersecting, subsets of the oriented circle, modulo orientation preserving homeomorphisms. The subsets $c_k$ are called chords.  We usually specify a chord diagram by drawing it as in figure \ref{chord}. In all the pictures we assume the circle oriented counterclockwise.

\begin{figure}
\centerline{\relabelbox
\epsfysize 1.5cm
\epsfbox{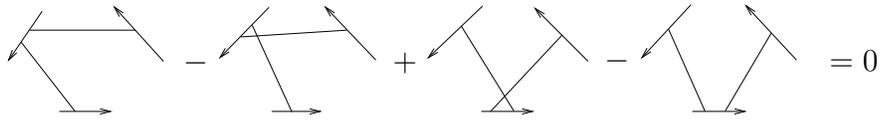}
\relabel {a}{$-$}
\relabel {b}{$+$}
\relabel {c}{$-$}
\relabel {d}{$=0$}
\endrelabelbox }
\caption{\label{4T} $4$ Term Relations.}
\end{figure}

For each $n \geq 2$, let $V_n$ be the free $\C$ vector space on the set of all
chord diagrams with $n$ chords. That is the set of formal finite linear
combinations $w=\sum_i \l_i w_i$, where $\l_i \in \C$ and $w_i$ is a chord diagram with $i$ chords for
any $i$. Consider the sub vector space $4T_n$ of $V_n$ which is the subspace
generated by all linear combinations of chord diagrams of the form displayed
in figure \ref{4T}. The $3$ intervals considered in the circle can appear in
an arbitrary order in $S^1$. Define for each $n \in \Nz=\{0,1,2,..\}$,the vector space $\A_n=V_n/4T_n$. We consider $\A_0=V_0$ and $\A_1 =V_1$.

\begin{figure}
\centerline{\relabelbox
\epsfysize 1.5cm
\epsfbox{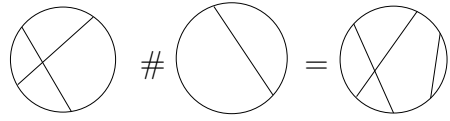}
\relabel {a}{$\#$}
\relabel {b}{$=$}
\endrelabelbox }
\caption{\label{connected} Connected Sum Product.}
\end{figure}

For any pair $m,n \in \Nz$, there exists a bilinear map $\#:\A_n \tn \A_m \to
A_{m+n}$, called the connected sum product. As its name says, it is performed
by doing the connected sum of chord diagrams as in figure \ref{connected}.

Obviously the  product is not well defined in $V_m \tn V_n$ since it depends on
the points in which we break the circles. The connected sum makes sense only in  $\A_m \tn
\A_n$, since we are considering the $4$-term relations. It is associative, commutative  and it has a unit: the chord diagram without any chord. For more details see \cite{BN}.

The vector space $\A_m \tn \A_m$ is mapped via the connected sum product to
$\A_{m+n}$. Therefore the direct sum $\A_\fin=\bigoplus _{n \in \Nz} \A_n$ has a
commutative and associative graded algebra structure. This permits us to
conclude that the vector space 
$$\A =\prod_{n\in \Nz} \A_n$$ 
has a structure of abelian algebra over the field of complex numbers. Call it the algebra of chord diagrams. The algebra $\A$ is the target space for the Kontsevich Universal Knot Invariant. 

\begin{figure}
\centerline{\relabelbox
\epsfysize 4cm
\epsfbox{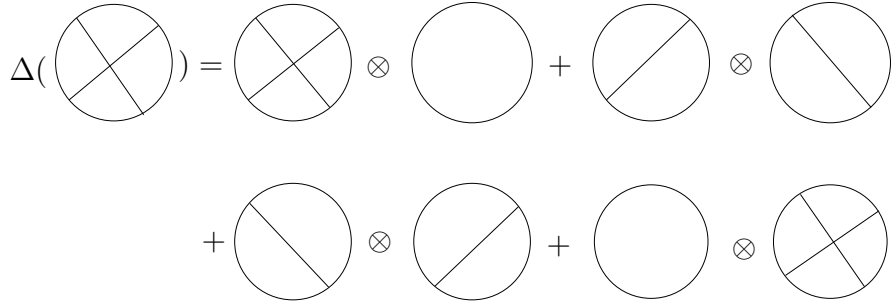}
\relabel {i}{$\Delta ($}
\relabel {a}{$)=$}
\relabel {b}{$\tn$}
\relabel {c}{$+$}
\relabel {d}{$\tn$}
\relabel {e}{$+$}
\relabel {f}{$\tn$}
\relabel {g}{$+$}
\relabel {h}{$\tn$}
\endrelabelbox }
\caption{\label{delta} Coproduct Maps.}
\end{figure}

There exist also coproduct maps $\D:\A_m \to \op_{k+l=m}A_k \tn \A_l$ which have the  form of figure \ref{delta} on chord diagrams.
They extend to a linear map $\D:\A \to \A \ctn \A$. Here $\A \ctn \A$ is the vector space 
$$\prod_{m \in \Nz} \bigoplus _{k+l=m} \A_k \tn \A_l.$$
Notice that $\A \tn \A$ is a proper sub vector space of $\A \ctn \A$.

An element $w \in \A$ is called group like if $\D(w) = w \ctn w$. That is, if writing $w =\sum_{n \in \Nz} w_n$ with $n \in \A_n, \forall n \in \Nz$ we have
$$\D(w_n) = \sum_{l+k=n} w_k \tn w_l.$$ 
For example, $\exp(\ominus)$ is a group like element. Here $\ominus$ is the
unique chord diagram with only one chord. This is a trivial consequence of the
fact $\D(\ominus)=\ominus \tn 1 +1 \tn \ominus$. We have put $1$ for the chord
diagram without chords.
\subsection{The Kontsevich Integral}
We skip the definition of the (framed) Kontsevich Integral $\cZ$, for which we refer for example to \cite{K}, \cite{LM} or \cite{SW}. See also \cite{BN,CV1} for the definition of the unframed version of the also called Kontsevich Universal Knot invariant. We take the normalisation of the Kontsevich Integral for which the value of the unknot is the wheels element $\W$ of \cite{NLT}. That is $\cZ(O)=\bf{Z}(\infty)$, cf \cite{BN} pp 447. This is a different normalisation of the one used in \cite{BN}. 
We now gather the properties of the Kontsevich integral which we are going to use in the sequel:

\begin {Theorem} \label{Kontsevich}
There exists a (oriented and framed) Knot invariant $K \mapsto \cZ(K)$, where $\cZ(K)$ is in the algebra $\A$ of chord diagrams. Given a framed knot $K$, $\cZ(K)$ satisfies:
\begin{enumerate}
\item  $\cZ(K)$ is grouplike, cf \cite{BN}
\item If $K^f$ is obtained from $\K$ by changing its framing by a factor of
  $1$ then $\cZ(K^f)=\cZ(K) \# \exp(\ominus)$, cf \cite{LM}. 
\item If $\K^*$ is the mirror image of $K$, and writing $\cZ(K)=\sum_{n \in \Nz} w_n$ with $\w_n \in \A_n, \forall n \in \Nz$ we have $\cZ(K^*)=\sum_{n \in \Nz}(-1)^n w_n$, cf \cite{CV1}.
\item If $K^-$ is the knot obtained from $K$ by reversing the orientation of it then $\cZ(K^-)=\sum_{n \in \Nz} S(w_n)$. Here $S:\A_n \to \A_n$ is the map that reverses the orientation of each chord diagram, cf \cite{CV1}.
\end{enumerate}
\end{Theorem}

Suppose we are given a family of linear maps (weights) $W_n :\A_n \to \C, n \in \Nz$. A knot invariant whose value on each knot is a formal power series with coefficients in $\C$ is called canonical if it has the form 
$$K \mapsto \sum_{n \in \Nz} W_n(w_n) h^n.$$
As usual we write $\cZ(K)=\sum_{n \in \Nz} w_n$ with  $w_n \in \A_n, \forall n \in \Nz$.

\subsection{Infinitesimal R-matrices}\label{irmatrices}

Let $\lg$ be a Lie algebra over the field $\C$. An infinitesimal R-matrix of
$\lg$ is a symmetric tensor $t \in \lg \tn \lg$ such that $[\D(X),t]=0,
\forall X \in \lg$. The commutator is taken in  $U(\lg) \tn U(\lg)$, where
$U(\lg)$ denotes the universal enveloping algebra of $\lg$. The map
$\D:U(\lg) \to U(\lg) \tn U(\lg)$ is the standard coproduct in
$U(\lg)$. It verifies $\D(X)=X\tn 1 +1\tn X$ if $ X \in \lg$.

Suppose we are given an infinitesimal R-matrix $t$. Write $t=\sum_i a_i \tn b_i$. We will then have:
$$\sum_{i,j} a_j a_i \tn b_i \tn b_j- a_i a_j \tn b_i \tn b_j+ a_i \tn a_j b_i \tn b_j- a_i \tn b_i a_j \tn b_j=0,$$
which resembles the $4T$ relations considered previously. Given a chord
diagram $w$ and an infinitesimal R-matrix $t=\sum_i a_i \tn b_i$ it is natural thus to construct an element $\f_t(w)$ of $U(\lg)$ in the following fashion: Start in an arbitrary
point of the circle and go around it in the direction of its
orientation. Order the chords of $w$ by the order with which you pass them as
in figure \ref{chords}. Each chord has thus an initial and an end
point. Then go around the circle again and write (from the right to the left)
$a_{i_k}$ or $b_{i_k}$ depending on whether you got to the initial or final
point of the $k^{th}$ chord. Then sum over all the $i_k$'s. For example for
the chord diagram of figure  \ref{chords} the element $\f_t(w)$ is:
$$\sum_{i_1,i_2,i_3} b_{i_2} b_{i_3} b_{i_1} a_{i_3} a_{i_2} a_{i_1}.$$
See \cite{K} or \cite{CV2} for more details. It is possible to prove that
$\f_t(w)$ is well defined as an element of $U(\lg)$, that is it does not depend
on the starting point in the circle. Moreover:

\begin{figure}
\begin{center}
\includegraphics{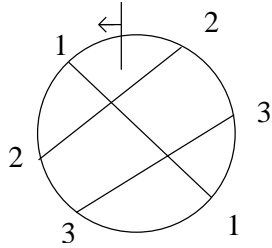}
\caption{\label{chords} Enumerating the Chords of a Chord Diagram.}
\end{center}
\end{figure}

\begin{Theorem} \label{central} Let $\lg$ be a Lie algebra and $t \in \lg \tn \lg$ be an infinitesimal R-matrix. The linear map $\f_t:V_n \to U(\lg)$ satisfies the $4T$ relations, therefore it descends to a linear map $\f_t:\A_n \to U(\lg)$. Moreover:
\begin{enumerate}
\item The image of $\f_t$ is contained in $\Ct(U(\lg))$, the centre of $U(\lg)$.
\item The degree of $\f_t(w)$ in $U(\lg)$ with respect to the natural
  filtration of $U(\lg)$ is not bigger than twice the number of chords of $w$. 

\item Given $w \in \A_m$ and $w' \in \A_n$ we have $\f_t(w \# w')=\f_t(w) \f_t(w')$

\item Consider the map $\f_{t,h}:\A \to \Ct(U(\lg))\h$ such that if
  $w=\sum_{n\in \Nz} w_n$ with $w_n \in \A_n$ for each $n \in \Nz$ we have 
$$\f_{t,h}=\sum_{n\in\Nz} \f_t(w_n)h^n.$$ Then $\f_{t,h}$ is a $\C$-algebra
morphism. 
\end{enumerate}

\end{Theorem}
Recall that the Kontsevich integral is a sum of the form $\cZ(K)=\sum_{n\in
  \Nz} w_n$ with $w_n \in \A_n, \forall n \in \Nz$. Therefore, given an
infinitesimal $R$ matrix $t$ in a Lie algebra $\lg$, we can obtain a knot invariant $\cZ_t$ which is defined as being
$$\cZ_t(K)=\f_{t,h}(\cZ(K))=\sum_n \f_t(w_n)h^n.$$
The target space of $\cZ_t$ is therefore the $\C$-algebra of formal power series over the centre of $U(\lg)$.

Suppose we are given a morphism $f:\Ct(U(\lg)) \to \C$. Then composing it with $\cZ(K)$ we obtain a canonical knot invariant $f \circ \cZ_t$. That is:
$$(f \circ \cZ_t)(K) =  \sum_n f(\f_t(w_n))h^n.$$

It is not difficult to examine the conditions whereby this kind of knot
invariants are unframed. Let $t=\sum_i a_i \tn b_i$ be an infinitesimal $R$
matrix in a Lie algebra. Define $C_t=\sum_i a_i b_i=-\f_t(\ominus)$. It is a
central element of the universal enveloping algebra of $\lg$. Call it the
quadratic central element associated with $t$. The infinitesimal $R$-matrix
$t$ can be recovered from $C_t$ by the formula
$$t=\frac{\D(C_t)-1\tn C_t-C_t \tn 1}{2}.$$
A morphism
$f:\Ct(U(\lg)) \to \C$ is
said to be $t$-unframed if $f(C_t)=0$. From theorem \ref{Kontsevich}, $2$ and
theorem \ref{central}, $3$ 
it is straightforward to conclude that:
\begin{Theorem}\label{unframed}
Let $\lg$ be a Lie algebra with an infinitesimal $R$-matrix $t$. Consider also
a morphism $f$ from the centre of $U(\lg)$ to $\C$. Then the knot invariant
$f\circ \cZ_t$ is unframed if and only if the morphism $f$ is $t$-unframed.
\end{Theorem} 
Notice the Kontsevich integral of each knot is invertible in $\A$. This is 
because the term $w_0 \in \A_0$ is the unit of $\A$.

\subsubsection{Constructing Infinitesimal $R$-matrices}\label{Constructing}

There exists a standard way to construct infinitesimal $R$-matrices in a Lie
algebra $\lg$. Suppose we are given a $\lg$-invariant, non degenerate,
symmetric bilinear form $<,>$ in $\lg$. Here $\lg$-invariance means that we
have  $<[X,Y],Z>+<Y,[X,Z]>=0,\forall X,Y,Z \in \lg$. If $\lg$ is semisimple the
Cartan-Killing form verifies the properties above. Take a basis
$\{X_i\}$ of $\lg$ and let  $\{X^i\}$ be the dual basis of $\lg ^*$ . Then it
is easy to show that for any $\l \in \C$ the tensor $t=\l \sum_i X_i \tn X^i$ is an infinitesimal $R$-matrix of
$\lg$. We are identifying $\lg^*$ with $\lg$ using the nondegenerate bilinear
form $<,>$.

Suppose $\lg$ is a semisimple Lie algebra and  let $t=\sum_i a_i \tn b_i$ be
an infinitesimal $R$-matrix in $\lg$. Let also $<,>$ denote the Cartan-Killing
form on $\lg$. Then the map $\lg \to \lg$ such that $X \mapsto \sum_i <X,a_i>
b_i$ is an intertwiner of $\lg$ with respect to its adjoint
representation. Therefore if $\lg$ is simple it is a multiple $\l$ of the
identity. This permits us to conclude that $t=\l X_i \tn X^i$. 

Let us now look at the case $\lg$ is semisimple. Then $\lg$ has a unique
decomposition of the form $\lg \cong \lg_1 \oplus ... \oplus \lg_n$, where each
$\lg_i$ is a simple Lie algebra. The Cartan-Killing form in each $\lg_i$
will yield an infinitesimal $R$-matrix $t_i$ in each $\lg_i$. Obviously each linear combination $t=\l_1 t_1+...+\l_n t_n$ is an infinitesimal $R$-matrix for
$\lg$. An argument similar to the one before proves that any infinitesimal
$R$-matrix in $\lg$ is of the form above. 

It should be said that in the case in which an infinitesimal $R$-matrix in a
Lie algebra $\lg$ comes from a non-degenerate, symmetric and $\lg$-invariant
bilinear form then our construction of central elements yields the same result
of \cite{BN}, cf \cite{CV2}.
\subsubsection{A Factorisation Theorem} 

Suppose the  Lie algebra $\lg \cong \lg_1 \op \lg_2$ is
the direct sum of two  Lie algebras. If $t_1$ and $t_2$ are infinitesimal
$R$-matrices in $\lg_1$ and $\lg_2$ then $t=t_1+t_2$ is also an infinitesimal
$R$-matrix in $\lg$. It is easy to prove that given a chord diagram $w$ we have the following
identity: cf (\cite{BN})

$$\f_t(w)=(\f_{t_1} \tn \f_{t_2})\D(w).$$ 
We are obviously considering the standard isomorphism $U(\lg) \cong U(\lg_1)
\tn U(\lg_2)$ such that $(X,Y) \mapsto X\tn 1 +1 \tn Y$ for $(X,Y) \in \lg$.

If we are given two algebra morphisms $f_i:\Ct(U(\lg_i)) \to \C,i=1,2$, then
$f=f_1 \tn f_1$ is an algebra morphism  $\Ct(U(\lg)) \cong
\Ct(U(\lg_i))\tn \Ct(U(\lg_i) \to \C$. It thus makes sense to consider the knot
invariant $f\circ \cZ_t$. It expresses in a simple form in terms of $f_i \circ
\cZ_{t_i}, i=1,2$. In fact: see \cite{melvin}

\begin{Theorem}\label{factorisation} 
Given any (oriented and framed) knot $K$ we have:
$$(f \circ\cZ_{t})(K)= (f_1 \circ
\cZ_{t_1})(K) \times (f_2 \circ \cZ_{t_2})(K), $$
as formal power series
\end{Theorem}
\begin{Proof}

Let $K$ be a knot, write $\cZ(K)=\sum_{n\in \Nz} w_n$ with $w_n \in \A_n,
\forall n \in \Nz$. We have:
\begin{align*}
(f \circ\cZ_{t})(K)&= \sum_{n\in \Nz} (f \circ \f_t) (w_n) h^n \\
&= \sum_{n\in \Nz} (f_1 \tn f_2)  \circ (\f_{t_1} \tn \f_{t_2}) (\D (w_n)) h^n \\
&= \sum_{n\in \Nz} \sum_{k+l=n} (f_1 \tn f_2)  \circ (\f_{t_1} \tn \f_{t_2})
(w_k \tn w_l) h^n \\
&= \sum_{n\in \Nz} \sum_{k+l=n}\left [(f_1   \circ \f_{t_1}) 
(w_k)\right ]\left [   (f_2   \circ \f_{t_2}) (w_l)\right ] h^{k+l} \\
&= (f_1 \circ \cZ_{t_1})(K)\times (f_2 \circ \cZ_{t_2})(K)
.\end{align*}  
\end{Proof}
This proof appears in  \cite{melvin}.         
\subsection{The Coloured Jones Polynomial}\label{Jones}

Let $\lg$ be a semisimple Lie algebra over $\C$. It is a well known result,
see for example \cite{VAR},
that any algebra morphism $\Ct(U(\lg)) \to \C$ is the central character of some
representation of $\lg$, which can be infinite dimensional. Recall that
$\Ct(U(\lg))$ stands for the centre of $U(\lg)$. To be more
precise, let $\lg$ be any Lie algebra and $\r$ a representation of $\lg$ in
the vector space $V$. Then $\r$ is said to admit a central character if every
element of $\Ct(U(\lg))$ acts on $V$ as a multiple of the identity. In this case there
exists an algebra morphism $\l_\r:\Ct(U(\lg)) \to \C$ such that $\r(a)(v)
=\l_\r(a) v, \forall a \in \Ct(U(\lg)), v \in V$. The algebra morphism $\l_\r$
is called the central character of the representation $\r$. 
In particular, if $\lg$ is a Lie algebra with an infinitesimal $R$-matrix $t$ then
given any representation $\r$ of $\lg$ with a central character, we can construct the knot invariant $(\l_\r \circ \cZ_t)$.

The Coloured Jones Polynomial is, up to normalisation,  a particular example of this
construction. Let $t$ be the infinitesimal $R$-matrix of $\slc$ corresponding
to the bilinear form in it which is minus the Cartan-Killing form. Consider for any
$\a \in \hN$ the representation $\ra$ of $\slc$ with spin $\a$, thus $\ra$
admits a central character which we denote by $\l_\a$. Given a framed knot $K$ Let $J^\a(K)$
denote the framed Coloured Jones Function of it. Notice we "colour"
the Jones polynomial with the spin of the representation, rather than with the dimension of it. The last one  is the usual convention. We have:
$$\frac{J^\a(K)}{2\a+1}=(\l_{\a}\circ \cZ_t)(K),\forall \a \in \hN.$$
Write 
$$\frac{J^\a(K)}{2\a+1}=\sum_{n\in \Nz} J^\a_n(K) h^n.$$
It is a known result that given a knot $K$ then $J^a_n(K)$ is a polynomial in
$\a$ with degree at most $2n$, cf \cite{MM}, \cite{C}. This is a consequence
of the fact that the centre of $U(\slc)$ is generated by the Casimir element
of it, together with $2$ of Theorem \ref{central}. Therefore we can
write:

$$\frac{J^\a(K)}{2\a+1}=\sum_{n \in \Nz}
\sum_{k=0}^{2n}a^{(n)}_k(K) \a ^k h^n.$$ 
For any complex number $z$ it thus makes sense to consider the $z$-Coloured
Jones Function. That is:
 $$\frac{J^z(K)}{2z+1}=\sum_{n\in \Nz} P^n(K)(z) h^n.$$ 
This yields thus a knot invariant whose value in a knot is a formal power
series in two variables: 
$$K \mapsto \sum_{m,n \in \Nz} a^{(n)}_k(K) z ^k h^n,$$ 
with $a^{(n)}_k(K) =0$ for $k>2n$.  
It is an interesting task to investigate whether or not this kind of series defines an analytic function in two variables. As we mentioned in the introduction they have in general a zero radius of convergence, so this can only be made precise in a perturbation theory point of view, cf \cite{FM}. This relates
to the question of whether it is possible to define \textit{numerical} knot
invariants out of the infinite dimensional representations of the Lorentz
Group. Notice it is known that if
$\a$ is a half integer then:
$$\frac{J^\a(K)}{2\a+1}= \sum_{m \in \Nz} \left (\sum_{k=0}^{2n} a^{(n)}_k(K) \a ^k \right ) h^n,$$
defines an analytic function in $h$.

For the unknot $O$ the series  $J^z(O)/(2z+1)$ has a non zero radius of convergence at any point $z \in \C$. The proof is not very difficult for we can have an explicit expression for it. Define, for each $z \in \C$ the meromorphic function:

$$F_z(h)=\frac{1}{2z+1}\frac{\sinh((2z+1)h/2)}{\sinh(h/2)}.$$
Thus for each $\a \in \hN$ we have 

$$F_\a(h)=\frac{J^\a(O)}{2\a+1}=\sum_{n \in \Nz} J^\a_n(O)h^n.$$
Consider the expansion:
$$F_z(h)=\sum_{n\in \Nz} c(z)_n h^n.$$
It is not difficult to conclude that each $c(z)_n$ is a polynomial in $z$ for fixed $n$. Moreover  $c(\a)_n=J^\a_{n}(O), \forall \a \in \hN$. This implies 
$$\frac{J^z(O)}{2z+1}=F_z(h),$$
as power series in $h$. In particular the power series for the unknot are convergent. This means it makes sense to speak about the quantum dimension of the representations of spin $z$, which are going to be defined later. To be more precise we made sense of the quantum dimension of them divided by their dimension as vector spaces. But notice the dimension of a representation of spin $z$ with $z \notin \hN$ is infinite. For some more explicit  examples see \cite{FM}.

\subsubsection{A Representation Interpretation of the $z$-Coloured \\ Jones Polynomial} \label{representation}

We can give an interpretation of the $z$-Coloured Jones Polynomial in the
framework of central characters. To this end, define the following elements of $\slc$:

 $$H=\begin{pmatrix}1 & 0 \\ 0 & -1 \end{pmatrix},  E=\begin{pmatrix}0 & 1 \\ 0
   & 0 \end{pmatrix}, F=\begin{pmatrix}0 & 0 \\ 1 & 0 \end{pmatrix}.$$
Then the infinitesimal $R$-matrix which we are considering in $\slc$ expresses
in the form:
$$t=-\frac{1}{4}\left (E \tn F +F \tn E +\frac{H \tn H}{2} \right ).$$
Notice that $t$ is defined out of the inner product in $\slc$ which is minus the Cartan-Killing form. In particular, the Casimir element $C$ of $\slc$ is equal to $-C_t$, where $C_t$ is
the quadratic central element associated with $t$. Recall subsection
\ref{irmatrices}.

Given a half integer $\a$, the representation space $\Va$ of the
representation of spin $\a$ has a basis of the form
$\{v_0,...,v_{2\a}\}$. The action of the elements $E,F$ and $H$ of $\slc$ in
$\Va$ is:

$$H v_k=(k-\a) v_k,$$
$$E v_k=(2\a-k) v_{k+1},$$ and
$$F v_k=k v_{k-1}$$
For an arbitrary complex number $z \notin \hN$, it makes sense also to speak
about the representation $\rz$ of spin $z$. Consider  $\Vz$ as being the infinite
dimensional vector space which has the basis $\{v_{2z},v_{2z-1},v_{2z-2},...\}$. Then
the representation $\rz$ of spin $z$ can be defined in the form:

$$H v_k=(k-z) v_k;k=2z,2z-1,...$$
$$E v_k=(2z-k) v_{k+1};k=2z,2z-1,...$$
$$F v_k=k v_{k-1};k=2z,2z-1,...$$

The  representations of spin $z \notin \hN$  have a central character
$\l_z$, since it is easily proved that each intertwiner $\Vz \to \Vz$ must be a
multiple of the identity. But see \cite{VAR}, 4.10.2., namely they are the unique irreducible cyclic highest weight representations with maximal weight $z$, this relative to the usual Borel decomposition of $\slc$.
Consider, given $z \in \C$, the framed knot invariant
$(\l_z \circ \cZ_t)$. Where, if $\a$ is half integer, $\l_\a$ is the central character of the usual representation of spin $\a$. Given a framed knot $K$ it has the form:

$$(\l_z \circ \cZ_t)(K)=\sum_{n \in \Nz} R_n^z(K) h^n,$$
where, by definition:

$$R_n^z(K)=(\l_z \circ \f_t)(w_n)=\sum _{n \in \N_0} \l_z (\f_t(w_n))h^n,$$ 
for 
$$\cZ(K) =\sum_{n\in \Nz} w_n, w_n \in \A_n, \forall n \in \Nz.$$ 
Also $$\frac{J^{\a}(K)}{2z+1}=\sum_{n \in \Nz} R_n^\a(K) h^n,\forall \a \in \frac{1}{2} \mathbb{N}.$$
Suppose $w$ is a chord diagram with $n$ chords. Let us have a look at the dependence of $\l_z( \f_t(w))$ in $z$. It is not difficult to conclude that it is a polynomial in this variable of degree at most $2n$. This is a trivial consequence of the definition  of the central element $\f_t(w)$ as well as the kind of action of the terms appearing in the infinitesimal $R$-matrix $t$ in \Vz. See also \cite{VAR} or \cite{FM}. In particular if $K$ is a framed knot,  $R_n^z(K)$ is a polynomial in $z$. Since we also have  $R_n^\a(K)=J_n^\a(K), \forall \a \in \hN$, we can conclude:

$$\frac{J^{z}(K)}{2z+1}=(\l_z \circ \cZ)(K),$$which gives us an equivalent definition of the $z$-Coloured Jones Polynomial.

The central characters of the representations of imaginary spin are actually the infinitesimal characters, cf \cite{Kir}, of the unitary representations of $SL(2,\R)$ in the principal series, cf \cite{L}, with the same parameter. Notice however that the derived representation of them in $\mathfrak{sl}(2,\R) \tn _\R \C \cong \sl$ is not any of the representation of imaginary spin just defined. This is the point of view considered in \cite{FM}.

\section{Lorentz Group}
Let $\lg$ be a semisimple Lie Algebra. As proved by Drinfel'D in \cite{D}, there is a one to one correspondence between gauge equivalence classes of
quantised universal enveloping algebras $\H$ of $\lg$ over $\Ch$, cf \cite{K}, and infinitesimal $R$-matrices in
$\lg$. Let us be more explicit about this. It is implicit in the definition of
a quantised universal enveloping algebra $\H$ that there exists a $\C$-algebra
morphism $f:\H/h\H \to U(\lg)$. Having chosen such morphism, the canonical $2$-tensor of $\A$ is defined as 
$t=f((R_{21}R-1)/h)$. It is an infinitesimal $R$-matrix of $\A$. Here $R$
denotes the universal $R$-matrix of $\H$. If $\H$ quantises the pair $(\lg,r)$ where $r$ is a classical $r$-matrix in $\lg$, see \cite{CP}, then $t$ is the symmetrisation of $r$. Each quantised universal enveloping
algebra can be given a structure of ribbon quasi Hopf algebra, cf \cite{AC}, and therefore there is a knot invariant attached each finite
dimensional representation of it, or what is the same, of $\lg$. These knot invariants take their values in the ring of formal power series over $\C$. If the
representation used is finite dimensional and irreducible then it has a central
character. In particular the framework of last section can be applied, using for example the infinitesimal $R$ matrix $t$ which is the canonical $2$ tensor of $\H$. It is a deep result that with these choices the two approaches for knot invariants are the same, up to division by the dimension of the representation considered. To be more precise we need also to change the sign of the infinitesimal $R$-matrix $t$, cf \cite{K}. 

In the case in which we consider a $q$-deformation $\A$ of the universal enveloping 
algebra of a Lie algebra $\lg$, then no such classification of gauge equivalence classes of quantised universal enveloping algebras exists. But
sometimes it is possible to make sense of the formula for $t$. This is  because we have a $q$-parametrised family of \textit{braided} Hopf algebras that tends to the universal enveloping algebra of $\lg$ as $q$ goes to $1$, or alternatively because $\A$ quantises the pair $(\lg,r)$ where $r$ is an $r$-matrix in $\lg$.

As mentioned in the introduction, despite the fact that the $q$-Drinfeld-Jimbo quantised universal enveloping
algebras $U_q(\lg)$ of semisimple Lie algebras are not ribbon Hopf algebras, their
category of finite dimensional representations is a ribbon category. That is
they have formal $R$-matrices and ribbon elements, which  make sense
when acting in their finite dimensional representations. The target space for
the knot invariants in this context is the complex plane. These
\textit{numerical} knot invariants can be obtained, apart from rescaling, by
summing the powers series which appear in context of $h$-adic Drinfeld-Jimbo algebras. In other words by summing the power series that come out of the approach making use the Kontsevich Integral and using the infinitesimal $R$-matrix which is the heuristic canonical $2$-tensor $t$ of $U_q(\lg)$.  

Let us pass now to the Quantum Lorentz Group $\cD$ as defined in
\cite{BR1} and \cite{BR2}. It is a quantum group depending on a parameter $q
\in (0,1)$. As said in the introduction, we wish to analyse
the question of whether or not there exists a knot theory attached to the
infinite dimensional representations of it. The situation is more or less the
same as the case of $q$-Drinfeld-Jimbo algebras. Namely we have an heuristic
$R$-matrix which comes from its structure of a quantum double as well as a
heuristic ribbon element. It is possible to describe how they act in the
unitary representations of $\cD$. The situation is simpler if  the minimal spin of the representation is equal to zero,  in which case the representation is said to be balanced. Representations of this kind are called simple in \cite{NR}.  In this context, the ribbon element acts as the identity and therefore the knot invariants obtained will be unframed. These invariants express out of an infinite sum as we will see in section \ref{Buf}.

One natural thing to do would be analysing whether the "derivatives" of these
sums define or not Vassiliev invariants, or  whether is possible to make sense
of them, in the framework of
Kontsevich Universal Invariant. It is not difficult to find an expression
for the heuristic canonical $2$ tensor of the quantum Lorentz Group. Also the
unitary representations of the quantum Lorentz Group in the principal and
complementary series have a classical counterpart. They are infinite
dimensional representations of the Lie algebra of the Lorentz Group which
admit a central character and therefore the framework of the last section can
be used. This is the program we wish to consider now.

\subsection{The Lorentz Algebra}\label{lorentz}

Consider the complex Lie group $SL(2,\C)$. Its Lie algebra
$\mathfrak{sl}(2,\C)$ is a complex Lie algebra of dimension $3$. A basis of
$\slc$ is $\{\s_X,\s_Y,\s_Z\}$ where

$$\s_X=\frac {1}{2}\begin{pmatrix} i & 0 \\ 0 & -i \end{pmatrix},\s_Y=\frac
{1}{2}\begin{pmatrix} 0 & i \\ i & 0 \end{pmatrix},\s_Z=\frac
{1}{2}\begin{pmatrix} 0 & -1 \\ 1 & 0 \end{pmatrix}.$$
The commutation relations are:
$$[\s_X,\s_Y]=\s_Z,\quad[\s_Y,\s_Z]=\s_X,\quad[\s_Z,\s_X]=\s_Y.$$
We can also consider a different basis $\{H_+,H_-,H_3\}$, where
$$H_+=i\s_X-\s_Y,\quad H_-=i\s_X+\s_Y,\quad H_3=i\s_Z,$$
the new commutation relations being:
$$[H_+,H_3]=-H_+,\quad[H_-,H_3]=H_-,\quad[H_+,H_-]=2H_3.$$
Restricting 
the ground field with which we are working to $\R$, we obtain the $6$
dimensional real Lie algebra $\slc_\R$, the realification of $\slc$. It is isomorphic with the Lie algebra of the Lorentz Group.

\begin{Definition}
The Lorentz Lie Llgebra $L$ is defined as being the complex Lie algebra
which is the complexification of $\slc_\R$. That is $L=\slc_\R\tn_\R \C$. 
It is therefore a complex Lie algebra of dimension $6$. The Lorentz Algebra 
is the complex algebra $U(L)$ which is the universal enveloping algebra of 
the
complex Lie algebra $L$.
\end{Definition}

The set $\{\s_X,B_X=-i\s_X,\s_Y,B_Y=-i\s_Y,\s_Z,B_Z=-i\s_Z\}$ is a real 
basis
of $\slc_\R$, and
thus a complex basis of $L$. The commutation relations are:
\begin{align*}
[\s_X,\s_Y]&=\s_Z,       &[\s_Y,\s_Z]&=\s_X,           &[\s_Z,\s_X]&=\s_Y,\\
[\s_Z,B_X]&=B_Y,         &[\s_Y,B_X]&=-B_Z,            &[\s_X,B_X]&=0,\\
[\s_Z,B_Y]&=-B_X,       &[\s_Y,B_Y]&=0,              &[\s_X,B_Y]&=B_Z,\\
[\s_Z,B_Z]&=0,          &[\s_Y,B_Z]&=B_X,            &[\s_X,B_Z]&=-B_Y,\\
[B_X,B_Y]&=-\s_Z,        &[B_Y,B_Z]&=-\s_X,            &[B_Z,B_X]&=-\s_Y.
\end{align*}

We can also consider the basis $\{H_+,H_-,H_3,F_+,F_-,F_3\}$ of $L$, where :

$$H_+=i\s_X-\s_Y,\quad H_-=i\s_X+\s_Y,\quad H_3=i \s_Z,$$
$$F_+=iB_X-B_Y,\quad F_-=iB_X+B_Y,\quad F_3=i B_Z.$$

The new commutation relations being:
$$[H_+,H_3]=-H_+, [H_-,H_3]=H_-,[H_+,H_-]=2H_3,$$
$$[F_+,H_+]=[H_-,F_-]=[H_3,F_3]=0,$$
$$[H_+,F_3]=-F_+, [H_-,F_3]=F_-,$$
$$[H_+,F_-]=-[H_-,F_+]=2F_3,$$
$$[F_+,H_3]=-F_+,[F_-,H_3]=F_-,$$
$$[F_+,F_3]=H_+, [F_-,F_3]=-H_-,[F_+,F_-]=-2H_3.$$

The following simple theorem will be one of the most important in our
discussion.
\begin{Theorem}

There exists one (only) isomorphism of complex Lie algebras  $\t: \slc \oplus \slc \to L \cong
\slc_\R
\tn_\R\C $ such that:
$$\s_X \op 0  \mapsto \frac{\s_X-i  \s_X \tn i}{2}=\frac {\s_X+iB_X}{2},$$
$$ 0\op \s_X  \mapsto \frac{\s_X+i  \s_X \tn i}{2}=\frac {\s_X-iB_X}{2},$$
$$\s_Y \op 0  \mapsto \frac{\s_Y-i  \s_Y \tn i}{2}=\frac {\s_Y+iB_Y}{2},$$
$$ 0\op \s_Y  \mapsto \frac{\s_Y+i  \s_Y \tn i}{2}=\frac {\s_Y-iB_Y}{2},$$
$$\s_Z \op 0  \mapsto \frac{\s_Z-i  \s_Z \tn i}{2}=\frac {\s_Z+iB_Z}{2},$$
$$ 0\op \s_Z  \mapsto \frac{\s_Z+i  \s_Z \tn i}{2}=\frac {\s_Z-iB_Z}{2}.$$
And thus we have also a Hopf algebra  isomorphism $$\t:U(\slc) \tn U(\slc) 
\to U(L).$$.
\end{Theorem}
\begin{Proof} Easy calculations\end{Proof}

Given $X \in \slc$, define $X^l=\t(X \op 0), X^r=\t(0 \op X)$.
And analogously for $X \in \Uslc$.
We have:
$$H_+^l=\frac{H_++iF_+}{2},\quad H_-^l=\frac{H_-+iF_-}{2}, \quad
H_3^l=\frac{H_3+iF_3}{2},$$
$$H_+^r=\frac{H_+-iF_+}{2}, \quad H_-^r=\frac{H_--iF_-}{2}, \quad
H_3^r=\frac{H_3-iF_3}{2}.$$
Consider also $C^l=\t(C \tn 1)$ and $C^r =\t(1 \tn C)$, where $C$ is the Casimir element of $\slc$ defined in \ref{Jones}. The elements $C^l$ and $C^r$ are called Left and Right Casimirs and their explicit expression is:
\begin{multline*}
4 C^l= \frac{H_3^2-F_3^2}{2} +i\frac {H_3F_3}{2}+i\frac{F_3H_3}{2} \\
+\frac{H_+H_-}{4}+i\frac{H_+F_-}{4}+i\frac{F_+H_-}{4}-\frac{F_+F_-}{4}\\
+\frac{H_-H_+}{4} + i\frac{H_-F_+}{4}+i\frac{F_-H+}{4}-\frac{F_-F_+}{4},
\end{multline*}

\begin{multline*}
4 C^r= \frac{H_3^2-F_3^2}{2} -i\frac {H_3F_3}{2}-i\frac{F_3H_3}{2} \\
+\frac{H_+H_-}{4}-i\frac{H_+F_-}{4}-i\frac{F_+H_-}{4}-\frac{F_+F_-}{4}\\
+\frac{H_-H_+}{4}- i\frac{H_-F_+}{4}-i\frac{F_-H+}{4}-\frac{F_-F_+}{4}.
\end{multline*}

We can also consider the left and right image under $\t \tn \t$ of the infinitesimal R-matrix of $\Uslc$. We take now $t \in \slc \tn \slc$ as being the infinitesimal $R$-matrix coming from the Cartan-Killing form. That is minus the one considered in \ref{Jones}. These  left and right infinitesimal 
R-matrices are:

\begin{multline*}
4 t^l= \frac{H_3\tn H_3}{2}-\frac{F_3\tn F_3}{2} +i\frac {H_3 \tn 
F_3}{2}+i\frac{F_3\tn H_3}{2} \\
+\frac{H_+\tn H_-}{4}+i\frac{H_+\tn F_-}{4}+i\frac{F_+\tn 
H_-}{4}-\frac{F_+\tn F_-}{4}\\
+\frac{H_-\tn H_+}{4} + i\frac{H_-\tn F_+}{4}+i\frac{F_-\tn 
H+}{4}-\frac{F_-\tn F_+}{4},
\end{multline*}

\begin{multline*}
4 t^r= \frac{H_3\tn H_3}{2}  -\frac{F_3 \tn F_3}{2} -i\frac {H_3 \tn 
F_3}{2}-i\frac{F_3 \tn H_3}{2}
\\+\frac{H_+\tn H_-}{4}-i\frac{H_+\tn F_-}{4}-i\frac{F_+\tn 
H_-}{4}-\frac{F_+\tn F_-}{4}\\
+\frac{H_-\tn H_+}{4}- i\frac{H_-\tn F_+}{4}-i\frac{F_ -\tn H+}{4}-\frac{F_ 
-\tn F_+}{4}.
\end{multline*}
Any linear combination  $at^l+bt^r$ of the left and right infinitesimal 
R-matrices is  an infinitesimal $R$ matrix for $L$. We wish to consider the combination $t_L=t^l-t^r$. That is

$$t_L=i \frac{1}{4}H_3 \tn F_3 + \frac{1}{4}iF_3 \tn H_3 +\frac{i}{8}H_-\tn 
F_+
+\frac{i}{8}F_-\tn H_+ +\frac{i}{8} H_+ \tn F_-+ \frac{i}{8} F_+ \tn H_-.$$
Notice another expression of it:
\begin{multline*} 
t_L =\\\frac{i}{8}\left (B_X \tn \s_X + \s_X \tn B_X+B_Y \tn \s_Y + \s_Y 
\tn B_Y+B_Z \tn \s_Z + \s_Z \tn B_Z \right ).\end{multline*}
The quadratic central element of $U(L)$ associated with $t_L$ is:

$$C_L=C_{t_L}= i \frac{H_3 F_3}{4} +i\frac{F_3 H_3}{4} + i\frac{H_+F_-}{8} +
i\frac{F_+H_-}{8}+i\frac{H_-F_+}{8} +i\frac{F_-H_+}{8}.$$
The reason why we consider this particular combination of the left and right infinitesimal $R$-matrices is because it corresponds to the heuristic canonical two tensor of the Quantum Lorentz Group considered in \cite{BR2}. Notice it is the symmetrisation of the classical $r$-matrix of $\mathfrak{sl}(2,\C)_\R$, see \cite{BNR} page 19. See also \cite{FM}. We shall see later (theorem \ref{equivalence}) that it is the right one. 
\subsubsection{The Irreducible Balanced Representations of the Lorentz 
Group} \label{representations}

Let us be given a complex number $p=|p|e^{i\theta},0 \leq \theta < 2 \pi $
different from zero. We define once for all the square root $\sqrt{p}$ of $p$ as
being
$\sqrt{|p|e^{i\theta}}=\sqrt{|p|}e^{i\frac{\theta}{2}}$. For $m \in \Z$ define the set 
$W_m=\{p \in \C : |p| \notin \N_{|m|+1}\}$, where, in
general, $\N_{m}=\{m,m+1,...\}$, for any $m \in \N$.  Consider 
the set
$\cP=
\{(m,p): m \in \Z, p \in W_m\}$. Define, for any $\a \in
\N$ and $ (m,p) \in \cP$:
$$C_\a(m,p)=\frac{i}{\a}\sqrt { \frac{(\a^2-p^2)(\a^2-m^2)}{4\a^2-1}},$$
$$B_\a(m,p)=\frac{ipm}{\a(\a+1)}.$$
Thus $C_\a(m,p) \neq 0, \a = |m|+1, |m|+2..., \forall p \in W_m$.

Consider the complex vector space $$V(m)=\bigoplus _{\a \in \N_{|m|}}\Va$$where $\Va$ denotes the representation space of the representation of $\slc$ of spin $\a$. The set $\{\va_i, i=-\a,-\a+1,...,\a; \a \in \N_{|m|}\}$ is 
a basis of $V(m)$. Consider the inner product in $V(m)$ that has the basis\ 
above as an orthonormal basis. Define also $\bar{V}(m)$ as being the Hilbert 
space which is the completion of $V(m)$.

Given $(m,p) \in \cP$, consider the following linear operators acting on 
$V(m)$:

$$H_3 \va_k=k\va_k, k=-\a,-\a+1,...,\a, \a \in \N_{|m|}$$
$$H_- \va_k=\sqrt{(\a+k)(\a-k+1)}\va_{k-1},
k=-\a,-\a+1,...,\a$$
$$H_+ \va_k=\sqrt{(\a+k+1)(\a-k)}\va_{k+1}, k=-\a,-\a+1,...,\a,$$
\begin{multline*}
F_+\va_k=C_\a(m,p) \sqrt{(\a-k)(\a-k-1)}\quad \vam_{k+1}\\
-B_\a(m,p)\sqrt{(\a+k+1)(\a-k)}\va_{k+1}\\
+ C_{\a+1}(m,p) \sqrt{(\a+k+1)(\a+k+2)}\quad \vaM_{k+1},\\
k=-\a,-\a+1,...,\a, \a \in \N_{|m|},
\end{multline*}

\begin{multline*}
F_-\va_k=-C_\a(m,p) \sqrt{(\a+k)(\a+k-1)}\quad \vam_{k-1}\\
-B_\a(m,p)\sqrt{(\a-k+1)(\a+k)}\va_{k-1}\\
- C_{\a+1}(m,p) \sqrt{(\a+1)^2-k^2}\quad \vaM_{k-1},\\
k=-\a,-\a+1,...,\a, \a \in \N_{|m|},
\end{multline*}

\begin{multline*}
F_3\va_k=C_\a(m,p) \sqrt{\a^2- k^2}\quad \vam_{k}
-B_\a(m,p) k \va_{k}\\
- C_{\a+1}(m,p) \sqrt{(\a+1)^2-k^2}\quad \vaM_{k},\\
k=-\a,-\a+1,...,\a, \a \in \N_{|m|}.
\end{multline*}
Obviously we are considering $\va_k=0$ if $k >\a$ or $k <-\a$.
We have the following theorem, whose proof can be found in \cite{GMS}

\begin{Theorem} If $(m,p) \in \cP$, the operators $H_-,H_+,H_3,F_-,F_+,F_3$ 
define
an infinite dimensional representation of the Lorentz Algebra.
\end{Theorem}

Notice that the representations $(m,p)$ and $(-m,-p)$ are equivalent. This has a trivial proof.

Denote the representations above by $\{\r(m,p):(m,p) \in \cP\}$. One can prove with no difficulty that they have a central character, for any intertwiner $V(m) \to V(m)$ needs to send each space $\Va$ to itself and act on it has a multiple of the identity. Considering the action of $F_+$, for example, we conclude that the multiples are the same in each space $\Va$. Therefore

\begin{Theorem}
For any $(m,p) \in \cP$ the representation $\r(m,p)$ of $L$ has a central character $\l_{m,p}$.
\end{Theorem}

For any $(m,p) \in \cP$ the representation $\r(m,p)$ of $L$ can always be 
integrated to a representation $R(m,p)$ of the Lorentz Group in the completion $\bar{V}(m)$ of $V(m)$, or to be more precise of  
its connected component of the identity. The representation is unitary if 
and only if $p$ is purely imaginary, for any $m\in \Nz$,  in which case the 
representation is said to belong to the principal series, or if $m = 0 $ and 
$p \in [0,1)$ in  which case the representation is said to belong to the 
complementary series. The vector space $V(m)$ is contained in the space of smooth vectors, cf \cite{Kir}, of $V(m)$; thus  $\l_{m,p}$ is the infinitesimal character of $R(m,p)$, in the unitary case. This unifies the approach here with the approach in \cite{FM}.

 The parameter $m$ is called the minimal spin of the representation.
A representation is called balanced if the minimal spin of it is $0$. 
Balanced representations depend therefore on a parameter $p \in W_0$. Denote 
them by $\{\r_p,p \in W_0\}$. Two balanced representations $\r_p$ and $\r_q$ of $L$ are equivalent if and only if $p=q$ or $p=-q$. These representations were used in \cite{BC} for the construction of a spin foam model for Quantum Gravity. The extension of that work for their quantised counterpart was dealt with in \cite{NR}.

Since the representations $\{\r(m,p): (m,p) \in \cP\}$ have a central character, the 
left and right Casimirs defined in \ref{lorentz} act on $V(m)$ as multiples of the identity. This 
multiples are, as a function of $m$ and $p$ the following:
$\frac{p^2+2mp+m^2-1}{8}$ for $C^l$ and $\frac{p^2-2mp+m^2-1}{8}$ for $C^r$. Therefore:
\begin{Proposition} If the infinitesimal R-matrix on $U(L)$ is the tensor $t_L$ defined in \ref{lorentz}, then the central characters  $\{\l_p,p \in W_0\}$  of the balanced representations  are $t_L$-unframed. Recall the nomenclature introduced before theorem \ref{unframed}.
\end{Proposition}

This can obviously be proved without using the explicit expression of the 
action of the Casimir elements.

Notice also that we can consider the minimal spin of the representations considered to be also to be an half integer,
making the obvious change in the form of the representation. These kind of representations cannot be integrated to
representations of the Lorentz Group, even though they define 
representations
of $SL(2,\C)$. They are called two-valued representations of the Lorentz 
Group in \cite{GMS} .

\subsection{The Lorentz Knot Invariant}

Consider again the  infinitesimal $R$-matrix $t_L=t^l-t^r$ of the Lorentz Lie Algebra. We consider for each $(m,p)\in \cP$ the representation $\r(m,p)$ of $L$. It has a central character $\l_{m,p}$. We propose to consider the framed knot invariants $\{X(m,p):(m,p) \in \cP\}$, such that for any knot:

$$K \mapsto X(m,p,K)=(\l_{m,p} \circ \cZ_t)(K).$$ 
Recall the notation of \ref{irmatrices}.
Notice $X(m,p)=X(-m,-p)$ for the representations $\r_{m,p}$ and $\r_{-m,-p}$ are equivalent.

The value of $X(m,p)$ in a framed knot $K$ is therefore a formal power series with coefficients in $\C$. It is a difficult task to analyse the analytic properties of such power series. We expect they will be perturbation series for some numerical knot invariants that can be defined, cf \cite{FM}.

As we have seen, if $m=0$, that is in the case of balanced representations, the central character $\l_p$ is $t_L$-unframed. This is also the case for $p=0$. Notice we have an explicit expression for the action of the left and right Casimir elements of $L$. Therefore

\begin{Theorem} \label{UNFRAMED}
The knot invariant $X(m,p)$ with $(m,p) \in \cP$ is unframed if and only if $m=0$ or $p=0$.
\end{Theorem}

Obviously, for different combinations of the left and right infinitesimal $R$-matrices, the representations which have unframed central characters with respect to it are different. This gives us a way to define an unframed knot invariant out of any $(m,p) \in \cP$. But notice this can be done without changing the infinitesimal $R$-matrix $t_L$ of $L$, since we know how the invariants behave with respect to framing, cf theorem \ref{Kontsevich}.

\subsubsection{Finite Dimensional Representations}\label{finrep}

Let us now analyse the knot invariants that come out of the finite dimensional
representations of the Lorentz Group. We are mainly interested in the
representations which are irreducible.

Since we have the isomorphism $U(L) \cong U(\slc) \tn U(\slc)$, the finite dimensional
irreducible representation of $U(L)$, or what is the same of $L$, are classified by a pair $(\a,\b)$ 
of half integers. That is each finite dimensional irreducible representation of $L$ is of
the form $\ra \tn \rb$ as a representation of $U(L) \cong U(\slc) \tn U(\slc)$. 
There is an alternative way to construct these finite dimensional
representations that shows their close relation with the infinite dimensional
representations, \cite{GMS}. Let us explain how the process goes. It is very similar to the $\sl$ case.

Consider $m=\a-\b$ and $p=\a+\b+1$. Notice
that now $C_\a(m,p)\neq 0$ if $\ \a \in |m|,|m|+1,...,p$, and $C_p(m,p)=0$. The
underlying vector space for the representation with spins $(\a,\b)$ is
$V(m,p)=\substack{|m|\\ V} \tn \substack{|m|+1\\V} \tn...\tn \substack {p-1\\V}$, and the form of it is given 
exactly by the same formulae of the infinite dimensional representations. The equivalence of the representations is a trivial consequence of the Clebsh-Gordan formula. This construction gives us a finite dimensional representation $\r(m,p)$ for each pair  $(m,p)$ with $m,p \in \Z/2$ and $p-|m| \in \N_1$. It makes also sense for $|p|-|m| \in \Z$, making the appropriate changes. As before we have the equivalence $\r(m,p) \cong \r(-m,-p)$. 

Since we completed the sets $W_m$ defined at the beginning of \ref{representations}, we have a  representation $\r(m,p)$ of the Lorentz Algebra for each 
pair $(m,p)$ with $m \in \Z/2$ and $p\in \C$. All them have a central character $\l_{m,p}$, since the new representations considered are finite dimensional and irreducible. The finite dimensional representations give us framed knot invariant $X_\fin(m,p)$ for each pair $m,p \in \Z/2$ with  $|p|-|m|\in \N_1$. This
invariant is independent of the framing if and only if $m=0$, that is if
$\a=\b$.

Consider now the algebra morphisms $\l_{m,p}\circ \f_{t_L}:\A \to \C$, where $m \in \Z$ and $p \in \C$. The argument is now similar to the one in \ref{representation}. If we look at the expression of the representations $\r_{m,p}$, it is easy to conclude that given any chord diagram $w$ with $n$
 chords, the evaluation of $\l_{m,p}\circ \f_{t_L}(w)$  is for a fixed $m$ a polynomial in $p$ of degree at most $2n$. Notice that any factor of the form $C_\a(m,p)$ appears in the expression for $\l_{m,p}\circ \f_{t_L}(w)$ an even number of times. For the case of balanced representations, that is $m=0$, we can also prove that it is a polynomial in $p^2$. Also the value of the polynomials in $p=1$ is zero if $n>0$ for the pair with $m=0$ and $p=1$ yields the trivial one dimensional representation of $L$. We have proved:

\begin{Theorem}\label{finite}
Consider the framed knot invariants $\{X(m,p), m\in \Nz, p \in \C\}$. If we fix $m\in\Nz$ then the term of order $n$ in the expansion of $X(m,p,K)$ as a power series is  polynomial of degree at most  $2n$ in $p$. Here $K$ is any framed knot. If $m=0$ then only the even terms of it are non zero. Moreover the polynomials attain zero at $p=1$ for $n>0$. 
\end{Theorem}

Therefore, if we know the value of $X(m,p,K)$ for the finite dimensional
representations, that is if $|p|-|m| \in \N$ we can determine it for any value of the parameter $p$. This is similar to the $\slc$ case.

\subsubsection{Relation With the Coloured Jones Polynomial}
The relation between the Lorentz knot invariants that
come out from finite dimensional and infinite dimensional representations remarked after theorem \ref{finite} gives
us a way to relate the Coloured Jones Polynomial with the Lorentz Group
invariants. In fact:

\begin{Theorem} \label{JonesRelation}
Let $K$ be some oriented framed knot, $K^*$ its mirror
  image. Then for any $z,w \in \C$ with $z-w \in \Z$ we have:

$$\frac{J^z(K^*)}{2z+1}\times \frac{J^w(K)}{2w+1}=X(z-w,z+w+1,K),$$
as formal power series over $\C$. 
\end{Theorem}

\begin{Proof} For any  $m\in \Z/2$ and $x \in \C$, let $z(x,m)=m+x$ and $w(x,m)=-m+x$. Thus each pair $(z,w) \in \C^2$ with $z-w \in \Z$ is of the form $(z(x,m),w(x,m))$ for some $m$ and $x$. Fix $m\in \Z/2$. We want to prove:

$$\frac{J^{(m+x)}(K^*)}{2m+2x+1}\times \frac{J^{(-m+x)}(K)}{-2m+2x+1}=X(m,2x+1,K), \forall x \in \C$$
Each term of the formal power series at both sides of the equality is a
polynomial in $x$, thus we only need to prove that the equality is true if both $x-m$ and $x+m$ are half integers. That is if $x-m,x+m \in \hN$.

 Let $t$ the infinitesimal $R$ matrix in $\slc$ coming out of the Cartan-Killing form. Notice it is minus the one considered in \ref{Jones}. Let $\a$ be a half integer. Recall that for a framed knot $K$ we have:

$$\frac{J^\a(K)}{2 \a+1}=(\l_\a \circ \cZ_{-t})(K).$$
Therefore by Theorem \ref{Kontsevich}, $3$:
$$\frac{J^\a(K^*)}{2\a+1}=(\l_\a \circ \cZ_{t})(K),$$
since $\f_t(w)=(-1)^n \f_{-t}(w)$ if $w$ is a chord diagram with $n$ chords.

Let $K$ be a framed knot and $x$ be such that $\a=x-m$ and $\b=x+m$ are  half integers.  
We have by theorem \ref{factorisation}:
\begin{align*}
\frac{J^\a(K^*)}{2\a+1}\times \frac{J^\b(K)}{2\b+1}&=(\l_\a \circ \cZ_{t})(K)\times (\l_\b \circ \cZ_{-t})(K)\\
&=(\l_\a \circ \cZ_{t^l})(K)\times (\l_\b \circ \cZ_{-t^r})(K)\\
&=\left ((\l_\a \tn \l_\b)\circ \cZ_{t_L} \right )(K).
\end{align*}
Recall $t_L=t^l-t^r$.

Now, $\l_\a \tn \l_\b$ is the central character of the representation $\r_\a \tn \r_\b$ of $U(L) \cong U(\slc) \tn U(\slc)$. As we have seen before, this representation is equivalent to $\r(\a-\b,\a+\b+1)=\r(m,2x+1)$. Thus their central characters are the same. This proves 

$$\left ((\l_\a \tn \l_\b)\circ \cZ_{t_L} \right )(K)=(\l_{m,2x+1}\circ \cZ_{t_L})(K)$$
if both $x-m$ and $x+m$ are half integers, and the proof is finished.
\end{Proof}

We have the following simple consequences.

\begin{Corollary}
Given a framed knot $K$, then the term of order $n$ in the power series of $X(m,z,K)$ is a polynomial in $m$ and $z$ 
\end{Corollary}

\begin{Corollary}
If $O$ is the unknot, then $X(m,p,O)$ is a convergent power series.
\end{Corollary}

\begin{Corollary} \label{mirror}
For balanced representations, that is if $m=0$, the invariant $X(0,p)$ does not distinguish a knot from its mirror image.
\end{Corollary}

\begin{Corollary}
The framed knot invariants $X(m,p)$ are unoriented.
\end{Corollary}
\section{The Approach with the Framework of Buffenoir and Roche}\label{Buf}
The aim of this section is to give a sketch of how the Buffenoir and Roche  description of the infinite dimensional unitary representations of the
Quantum Lorentz Group relates with our approach. The Quantum Lorentz Group was
originally defined by Woronowicz and Podle\'{s} in \cite{PW}. The classification
of the irreducible unitary representations of it appeared first in  \cite{P}.

For an expanded treatment of the issues considered in this section, we refer
the reader to \cite{PHD}.

\subsection{Representations of the Quantum Lorentz Group and $R$-Matrix}
We now follow \cite{BR1}. Other good references are \cite{BR2} and \cite{BNR}.
These references contain all the notation and conventions we use. The Quantum
Lorentz Group $\cD$ at a point $q \in (0,1)$ is defined as the quantum double
$\cD(\su, \pol)$. Notice that both $\su$ and $\pol^{\mathrm{cop}}$ are sub
Hopf algebras of $\cD$. The Quantum Lorentz Group  thus have a
\textit{formal} $R$-matrix coming from its quantum double structure. Even
though it is defined by an infinite sum, it is possible to describe its action
in any pair of infinite dimensional irreducible representations of $\cD$ in
the principal series. See \cite{BR1,BR2}  for a description of them. For the
dual counterpart of the theory, in other words for the theory of
corepresentations of the algebra of function in the Quantum Lorentz Group
$SL_q(2,\C)$ we refer to \cite{PuW}. 

Let us describe what the situation is in the case the two  representations are the same. Suppose also the minimal spin $m$ of them is zero.  Similarly with the classical case described above, representations $\r(p)$ of this kind  will be called balanced. They  depend a  parameter $p \in \C$. If $p \in i\R$ then the representations $\r(p)$ can be made unitary. Choosing  $p \in [0,i \frac{2\pi}{h}]$, where $q=e^{h/2}$, parametrises all the unitary representations in the principal series which have minimal spin equal to zero. These last ones are called simple representations in \cite{NR}.

Similarly with the classical case, the underlying vector space for the balanced representations $\r(p), p\in \C$ of the Quantum Lorentz Group is $$V=V(p)=\bigoplus_{\a \in \Nz} \Va,$$
where $$\ra :\su \to L(\Va)$$
is the irreducible representation of $U_q(\mathfrak{su}(2))$ with spin $\a$.
A basis of $\Va$ is thus given by the vectors 
\hbox{${\{\va_i,i=-\a,-\a+1,...,\a\}}$}. Any element $x$ of $\su$ acts in $V$ 
  in the fashion:
$$\prod_{\a\in \Nz}\ra(x).$$
The group like element of the Lorentz Group is given by $G=q^{2J_z}$. The 
heuristic ribbon element of the Quantum Lorentz Group is easily proved to
act as the identity in the balanced representations. See \cite{NR}.

Define, given half integers $A,B,C$ and $D$, the complex numbers: 
\begin{equation}\label{DefLambda}
\La{A}{B}{C}{D}(\r)=\sum_{\s}\CGR{A}{C}{B}{0}{\s}{-\s} 
q^{2\s\r}\CGL{B}{C}{D}{-\s}{\s}{0}.
\end{equation}
For the correct definition of the phases of the Clebsch-Gordan coefficients  see \cite{BR2}. We display their explicit expression later. The formal universal $R$-matrix of the Quantum Lorentz Group is: (see \cite{BR1})
$$\mathcal{R}=\sum_{\substack{\a \in \hN\\-\a \leq i_a,j_a \leq \a}}\Xai\tn \gaj,$$
its inverse being:
$$\mathcal{R}^{-1}=\sum_{\substack{\a \in \hN\\-\a \leq i_a,j_a \leq \a}}\Xai\tn S^{-1}(\gaj).$$
This  antipode $S$ is the one of $\pol^{\mathrm{cop}} \subset \cD$, which is the inverse of the one in $\pol$, thus
 $$S^{-1}(\gaj)=q^{-i_\a+j_\a}(-1)^{-j_\a+i_\a}\ga^{-i_\a}_{-j_\a}.$$
See \cite{BR2}, equation $(25)$.
The action of $\gai$ in the space $V(p)$ is given by
\begin{equation}\label{gaction}
\gai \vbi=\frac{\F_\b}{\F_\g}\sum_{D,\g,x \hN}\sum_{-\g\leq \ig\leq \g}\sum_{-D\leq x\leq D}
\vgi\CGL{\g}{\a}{D}{i_\g}{i_\a}{x}\CGR{D}{\a}{\b}{x}{j_\a}{i_\b} 
\La{\g}{\a}{C}{\b}.
\end{equation}
Note it is a finite sum. The constants $\F_\a$ are defined in \cite {BR2}, proposition $1$. They will
not be used  directly. In fact their values are (almost) arbitrary and they
only appear to ensure that the representations $\r(p)$ are unitary for $p \in
i\R$, and the natural inner product in $V$. Their appearance does not  change the representation itself, therefore
not affecting the calculations of knot invariants.

The coefficients  $\La{A}{B}{C}{D}(p)$  are originally defined in \cite{BR2}
from an analytic continuation of $6j$-symbols, and at the end  proved to
coincide with (\ref{DefLambda}). One can show directly that (\ref{gaction})
does define a representation of the Quantum Lorentz Group, for any $p\in \C$
since equation $(76)$ of \cite{BR2} holds. See \cite{PHD}.

In some particular cases, equation (\ref{gaction}) simplifies to:
\begin{equation}\label{simplified1}
  \gai \vo=\frac{\F_0}{\F_\g}\sum_{\g,\ig}\CGL{\g}{\a}{\a}{i_\g}{\ia}{\ja}
  \La{\g}{\a}{\a}{0} \vgi,
\end{equation}
and
\begin{equation}\label{simplified2}
\left< \voo, \gai \vbi\right>=\frac{\F_\b}{\F_0}\CGR{\a}{\a}{\b}{\ia}{\ja}{\ib}
  \La{0}{\a}{\a}{\b}. 
\end{equation}
All these formulae are consequences of  well known symmetries of Clebsch-Gordan Coefficients listed for example in \cite{BR2}. With them we can also prove $\La{\a}{0}{\a}{\a}=1$, from which follows:
\begin{equation}\label{actiong0}
\go \vai=\vai.
\end{equation}

The elements $$\Xai\in \pol^*,\a \in \hN,i_\a={-\a,...,\a}$$ act simply as matrix
elements, that is:
\begin{equation}\label{xaction}
\Xai \vbi=\d(\a,\b)\d(i_a,i_\b)\vaj.
\end{equation}
Notice $\su$ is naturally embedded in $\pol^*$. Moreover any finite
dimensional irreducible representation of $\su$ induces one of $\sp\left
  \{\Xai\right \}\subset \pol^*$ which
has exactly this form, cf \cite{PW}, theorem $5.1$.

The action of the group like element $G$ is
\begin{equation}\label{Gaction}
G\vai=q^{2i_a}\vai.
\end{equation}
It is  easy to compute how $\mathcal{R}$ acts, namely:

$$\mathcal{R} \left (\vai \tn \vbi\right)=\sum_{D,x,\g,i_\g,j_\a}\CGL{\g}{\a}{D}{i_\g}{j_\a}{x}
\CGR{D}{\a}{\b}{x}{i_\a}{i_\b}\frac{\F_\b}{\F_\g}\La{\g}{\a}{D}{\b}\left ({\vaj}
\tn \vgi\right). $$
See \cite{BR1}, proposition $13$. The domain of the sum is the obvious one. The action of $\mathcal{R}$ in $V \tn V$ is thus well defined. Note  we are
considering the algebraic, rather than topological,  tensor product. Moreover  $\mathcal{R}$  defines a braid group 
representation. Denote it by $b \in B(n) \mapsto R_b \in L(V^{\tn})$. Here
$B(n)$ denotes the n-strand braid group and $L(V^{\tn n})$ the vector space of  linear maps $V^{\tn n} \to V^{\tn n}$. Notice that the braiding operators $R_b$ extend to unitary operators if  $p \in i\R$ since $\mathcal{R}^{*\tn *}=\mathcal{R}^{-1}$, where $*$ is the star structure on the Quantum Lorentz Group, see \cite{BNR}, since $\r(p)$ is unitary in this case.
\subsection{Associated Knot Invariants}
We now use the framework just introduced to define quantum lorentzian knot invariants. As we will see they relate to our approach before.
\subsubsection{Some Heuristics}
Let $q \in (0,1)$ and $p \in \C$. Suppose   we are given a braid $b$ with $n+1$ strands. There is  attached to it a map $R_b:V^{\tn (n+1)} \to V ^{\tn (n+1)}$. Consider the
map $A_b=(\id \tn G\tn...\tn G)R_b$. Suppose the closure of the braid $b$ is a
knot. If the representations we are
considering were finite dimensional, then the partial trace $T^1(A_b):V \to
V$ of $A_b$ over the last $n$ variables would be an intertwiner and thus a
multiple of the identity,  since the representations which we are
considering are irreducible. Moreover \textit{this} multiple of the identity
would be a knot invariant, which  would have the form:
\begin{multline}\label{FormalRT1}
b \mapsto S_b(q,p)=\\
\sum_{\substack{\a_1,...,a_{n} \in \hN \\-\a_k\leq i_{\a_k} \leq \a_k,k=1,...,n}}\left <\voo \tn \vaa^{i_{\a_1}} \tn... \tn
\van^{i_{\a_n}},A_b\left (\vo \tn \vaa_{i_{\a_1}} 
\tn... \tn \van_{i_{\a_n}}\right)\right >.
\end{multline}
 Even though the sums above may be not convergent, the assignment of one sum
of this kind to a braid whose closure is a knot is not ambiguous.  In fact
suppose $b$ has $m+1$ strands and $n$ crossings. We can always express this
sum in a more suggestive way, namely as:
\begin{equation}\label{FormalRT}
S_b(q,p)=\sum_{\substack{a_1,...,\a_n\\-\a_k \leq i_k,j_k\leq \a_k,k=1,...,n}} \left
  <\voo,\prod_{l=1}^{2n+m} T(\aa,\ii,\jj,l)\vo \right >,
\end{equation}
where if $\aa=(\a_1,...,\a_ n)$, $\ii=(i_1,...,i_n)$ and  $\jj=(j_1,...,j_n)$, then
$ T(\aa,\ii,\jj,l)$ can be either a term of the form  $\gak^i_j$ or $\Xak^i_j$, for some $k\in \{1,..,n\}$ or $G$; and moreover for any $k$ there exists an $l$ such that $ T(\aa,\ii,\jj,l)$ is a $\Xak^i_j$, and the same for $\gak^i_j$.  The two examples below should clarify what we mean. We obviously need to suppose that the closure of $b$ is a 
knot for this to hold. Notice that the transition from (\ref{FormalRT1}) to
(\ref{FormalRT}) is totally clear if the representations are finite
dimensional. We take (\ref{FormalRT}) as the definition of $S_b(q,p)$ if $b$ is a
braid whose closure is a knot.

 Let us look to the sums above in
a bit more detail. We consider the Left and Right Handed Trefoil
knots displayed in figure \ref{trefoil}. Call the two braids we have chosen
to represent them $T_+$ and $T_-$. The sum for the Right Handed Trefoil Knot is:

\begin{figure}
\centerline{\relabelbox
\epsfysize 3.5cm
\epsfbox{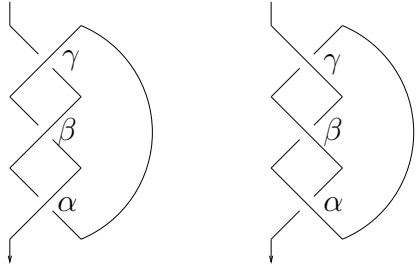}
\relabel {a}{$\g$}
\relabel {b}{$\b$}
\relabel {c}{$\a$}
\relabel {d}{$\a$}
\relabel {f}{$\b$}
\relabel {g}{$\g$}
\endrelabelbox}
\caption{\label{trefoil} Right and Left Handed
Trefoil Knots.}
\end{figure}

$$S_{T_+}(q,p)=\sum_{\substack{\a,\b,\g \in \hN \\ \a\leq i_\a,j_\a \leq \a \\\b\leq i_\b,j_\b \leq \b \\ \g\leq i_\g,j_\g \leq \g }}\left <\voo,\ggj\Xbi\gaj G \Xgi \gbj \Xai \vo\right>.$$
Whereas for the Left Handed Trefoil is:
$$S_{T_-}(q,p)=\sum_{\substack{\a,\b,\g \in \hN \\ \a\leq i_\a,j_\a \leq \a \\\b\leq i_\b,j_\b \leq \b \\ \g\leq i_\g,j_\g \leq \g }} \left <\voo,\Xgi S^{-1}(\gbj)\Xai G S^{-1}(\ggj) \Xbi S^{-1}(\gaj)
\vo\right>.$$
Many of the terms will be zero in the expressions above. Let us look 
at $S_{T_-}$. We only want the $0 \to 0$ matrix element, and
$\left <\voo,\Xai v \right >= \d(\a,0)\left <\voo,v\right>$. Thus we can make $\g=0$,  and then note that
$\go$ acts as the identity. We obtain: (we skip unnecessary indices)
\begin{multline*}
S_{T_-}(q,p)\\ =\sum_{\a,\b \in \hN}q^{-i_\a +j_\a-i_\b+j_\b}(-1)^ {i_\a+i_\b-j_\a-j_\b}\left <\gb^{-\ib}_{-\jb} \Xai G\Xbi \ga ^{-\ia}_{-j_a}\vo,\vo\right>.
\end{multline*}
From (\ref{xaction}) and (\ref{Gaction}) follows  $\a=\b$ and
$\ia=\jb$. By (\ref{simplified1}) and (\ref{simplified2}) we can conclude:
\begin{multline*}
S_{T_-}(q,p)=\sum_{\a \in \hN}\sum_{i_\b,\jb,j_\a=-\a}^\a q^{j_\a-i_\b +2\jb}(-1)^ { i_\b-\jb} \\ \CGR{\a}{\a}{\a}{-\ib}{-\jb}{\ja}
\CGL{\a}{\a}{\a}{\ib}{-j_\b}{-j_\a} \La{0}{\a}{\a}{\a} \La{\a}{\a}{\a}{0}.
\end{multline*}
Using the standard symmetries of the Clebsch-Gordan Coefficients,  we can express this as:
$$S_{T_-}(q,p)=\sum_{\a \in \Nz}\sum_{i_\b,\jb,j_\a=-\a}^\a q^{2\jb}
\CGL{\a}{\a}{\a}{-\ib}{-j_\a}{-\jb}\CGR{\a}{\a}{\a}{-\jb}{-\ib}{-j_\a}
\La{0}{\a}{\a}{\a} \La{\a}{\a}{\a}{0}.$$
Notice $\La{\a}{\a}{\a}{0}$ is zero unless $\a$ is integer. Therefore the final expression for the sum is:
\begin{equation*}
S_{T_-}(q,p)=\sum_{\a \in \Nz} d_\a \La{0}{\a}{\a}{\a} \La{\a}{\a}{\a}{0}.\end{equation*}
Here $d_a$ is the quantum dimension of the representation $\ra$. It equals $(q^{2\a+1}-q^{-2\a-1})/ (q-q^{-1})$. This last sum  is easily proved to be equal to $S_{T_+}$, therefore, if the sums do  define a knot invariant, they make no distinction between the Trefoil and its  mirror image. We would expect this from corollary \ref{mirror}. The calculations for other knot diagrams follow the same procedure, which can be given an obvious graphical calculus.

Notice that the series $S(T_-)$ seems to be divergent due to the presence of the $d_\a$ term in it. Therefore this sums do not seem to  define $\C$-valued knot invariants. This tells us the method of Borel re-summation sketched in \cite{FM} is perhaps more powerful.
\subsubsection{Finite Dimensional Representation}\label{qfinrep}
Let $Y(\a,\b,\g)=1$ if $\ra$ is in the decomposition of $\rb \tn \rg$ in term
of irreducible representations of $\su$ and zero otherwise, where $\a,\b,\g \in
\frac{1}{2}\Z$. Let also $Y(\a,\ia)=1$ if $\ia \in \{-\a,..,\a\}$ and zero otherwise. We have (see \cite{BR2}):
\begin{multline}\label{CGformula}
 \CGL{I}{J}{K}{m}{n}{p}=Y(I,m)Y(J,n)Y(K,p)\delta(m+n,p)Y(I,J,K)\\q^{m(p+1)+\frac{1}{2}((J(J+1)-I(I+1)-K(K+1)}e^{i\pi (I-m)}\\
\sqrt{\frac{[2K+1][I+J-K]![I-m]![J-n]![K-p]![K+p]!}  {[K+J-I]![I+K-J]![I+J+K+1]![I+m]![J+n]!}}\\
\times \sum_{\substack{V=0\\-J+K-m\leq V\leq I-m}}^{K-p}\frac{q^{V(K+p+1)} e^{i\pi V}[I+m+V]![J+K-m-V]!}{ [V]![K-p-V]![I-m-V]![J-K+m+V]!           }
.\end{multline}
Let $p\in \C$. We thus have an infinite dimensional representation $\r(p)$ of the Quantum Lorentz Group given by the constants $\La{A}{B}{C}{D}(p)$. Its representation space is by definition $V =V(p)=\bigoplus_{\a \in \N_0} \Va$. From equation (\ref{CGformula}), we can easily calculate the coefficients $\La{A}{B}{C}{D}$ if $B=1/2$. In fact
\begin{Lemma}\label{lfor}
Let $C\geq 0$ be an integer . We have
\begin{align*}
  \La{C}{1/2}{ C-1/2}{C}{(p)}&=\frac{q^{C}\left(  q^{p}+q^{-p}\right )}{q^{2C}+1},\\
\La{C}{1/2}{C+1/2}{C}{(p)}&=-\frac{q^{C+1}\left(q^{p}+q^{-p}\right
  )}{q^{2C+2}+1},\\
 \La{C}{1/2}{C+1/2}{C+1}{(p)}&=\frac{q^{2C+2}  q^{p}-q^{-p}}{q^{2C+2}+1},\\
 \La{C+1}{1/2}{C+1/2}{C}{(p)}&=\frac{q^{2C+2}  q^{-p}-q^{p}}{q^{2C+2}+1}
. \end{align*}
\end{Lemma}
Notice that all the other $\La{A}{B}{C}{D}(p)$ coefficients with $B=1/2$ are zero. See also  \cite{BR2}, proof of theorem $3$.

In particular, if $p \in \N$, then the representation $\r(p)$ has a finite dimensional subrepresentation $\r(p)_\fin$ in $V(p)_\fin=\Vo \tn \Vu \op...\op  \Vpm$. Compare with \ref{finrep}. Using Schur's lemma as in \cite{BR2}, proof of theorem $3$, one proves these representations are irreducible.

 Notice that $\vo\in V(p)_\fin$. Therefore, looking at (\ref{FormalRT}) we conclude:
\begin{Lemma}\label{truncate}
If $p\in \N$ and $b$ is a braid whose closure is a knot, then the infinite sum $S_b(q,p)$ truncates to a finite sum for any $q\in (0,1)$.
\end{Lemma}
As we referred in the introduction, the category of finite dimensional
representations of the Quantum Lorentz Group is (almost) ribbon equivalent to
the category of finite dimensional representations of $\su \tn_{R^{-1}} \su$
where $R$ is the $R$-matrix of $\su$. Let us explain what this means. We
follow \cite{BR2} and \cite{BNR} closely. The ribbon Hopf algebra $\su
\tn_{R^{-1}} \su$ is isomorphic with $\su \tn \su$ as an algebra, but has a
coalgebra structure of the form: 
$$\D(a \tn b)=R_{23}^{-1}a'\tn b' \tn a'' \tn b'' R_{23},$$  whereas the antipode is defined as:
$$S(a\tn b)=R_{21}S(a)\tn S(b) R_{21}^{-1}.$$
This Hopf algebra  has an $R$-matrix given by
$$\hat{R}=R_{14}^{(-)}R_{24}^{(-)}R_{13}^{(+)}R_{23}^{(+)},$$
where  $R^{(+)}=R$ and $R^{(-)}=R_{21}^{-1}$. The algebra $\su \tn_{R^{-1}} \su$ is a ribbon Hopf algebra with  group like element  $G\tn G$, where $G=q^{2J_z}$ is the group like element of $\su$. See \cite{BNR} page $20$.

The irreducible finite dimensional representations $\ra_w$ of $\su$ are
parametrised by an $\a\in \hN$ and an $w \in \{1,-1,i,-i\}$, the level of the
representation. See \cite{KS}, theorem $13$. The irreducible representations
$\ra$ of spin $\a$ are the ones for which $w=1$. They are the natural
quantisation of the representations of $SU(2)$ of spin $\a$. The action of the
$R$-matrix of $\su$ is, apriori, only defined on pairs of representations of
level $1$, or direct sums of them. Nevertheless, the category of finite
dimensional representations of $\su$ of this kind is a ribbon category.    
Let $\ra$ and $\rb$ be two finite dimensional irreducible representations of
$\su$ of level $1$ which will then generate a representation $\ra \tn
\rb$ of $\su \tn_{R^{-1}} \su$. The action of the $R$-matrix of  $\su
\tn_{R^{-1}} \su$ is well defined on pairs of representations of this
kind. The same is true for the group like element,   thus
we can define a framed  knot invariant  $I(\ra \tn \rb)$. Unpacking the expression of it yields immediately
\begin{Lemma}\label{mirrors}
For any framed knot $K$ we have
$$I(\ra \tn \rb)(K)=I(\ra)(K^*) I(\rb)(K)$$
where $K^*$ is the mirror image of $K$. Here $I(\ra)$ is the $\su$-framed knot invariant defined from $\ra$, in other words the Coloured Jones Polynomial, and the same for $I(\rb)$.
\end{Lemma}
There exists a Hopf algebra morphism $\psi:\cD \to \su \tn_{R^{-1}} \su$. It has the form:
$$\psi:(x,f)\mapsto \sum_{(x)(f)}x'(f''\tn \id)(R^{(+)})\tn x''(f'\tn \id)(R^{(-)})$$
The comultiplications are taken in $\su$ and $\pol$. This morphism naturally
extends to the elements $\Xai\in \pol^*$.

There also exists a morphism $s:\pol \to \pol$ such that $s(\gai)=(-1)^{2\a} \gai$. It
 extends to all the Quantum Lorentz Group provided we define its restriction
 to $\su$ (thus also to $\sp\left
  \{\Xai\right \}\subset \pol^*$ ) to be the identity.  
 A main result of  \cite{T}, namely theorem $5.4$,  is the following:
\begin{Theorem}
Let $\r$ be an irreducible finite dimensional representation of the 
Quantum Lorentz Group $\cD$ which has a structure of a $\pol$-crossed bimodule. In
our case this means that the representation $\r$ restricted to $\su$ is a
direct sum of representations $\ra_w$ with $\w=1$, which is what happens  for
the representations $\r(p)_\fin, p \in \N$, see \cite{T} proposition $5.1$ (These representations define
corepresentations of the algebra $S_qL(2,\C)$, thus representations of the
Quantum Lorentz Group in the sense of \cite{PW}). Then there exists $\a,\b \in
\hN$ such that either $\r=(\ra \tn \rb)\circ \psi$ or $\r= (\ra \tn \rb)\circ
\psi\circ s$.
\end{Theorem}
See \cite{BR2}, page $507$. Therefore if $p\in \N$  and $\a=(p-1)/2$, then
either $\r(p)=(\ra \tn \ra) \circ \psi$ or $\r(p)=(\ra \tn \ra_i) \circ \psi
\circ s$,
since the minimal spin of $\r(p)_\fin$ is zero (note that $\psi$ restricted to $\su$ is simply the coevaluation $\D$). With a bit more  work one
can actually prove that it is the second case that holds. See \cite{PHD}

 If we consider the action in finite dimensional representations of the form
 $\r=(\ra \tn \rb)\circ \psi$, then $\psi$ transforms the $R$-matrix of the
 Quantum Lorentz Group into the $R$-matrix of $\su \tn_{R^{-1}} \su$, and analogously for their inverses. This is
 an easy consequence of the fact $(\D \tn \id)(R)=R_{13}R_{23}$ and $(\id \tn
 \D)(R)=R_{13}R_{12}$. The same is true for the  balanced representations
 $\r(p)\cong(\ra \tn \ra_i) \circ \psi \circ s$ since given that $\Xai\tn \gaj$ acts as zero in $V(p)_\fin \tn V(p)_\fin$
 if $\a\in \N+\frac{1}{2}$, it follows that the actions of $(s \tn
 s)(\mathcal{R})$ and $\mathcal{R}$ in $V(p)_\fin\tn V(p)_\fin $ are the
 same. The detailed calculation appears in \cite{PHD}  
 The map $\psi$  preserves the group like elements since $\D(q^{2J_z})=q^{2J_z} \tn q^{2Jz}$. Therefore from lemma \ref{mirrors}  we obtain:
\begin{Proposition}\label{ThePoint}
Let $q\in (0,1)$ and $p\in \N$. Let also $\a=(p-1)/2$. Given a braid $b$, let
$K_b$ be the closure  of $b$ with an arbitrary framing and $K_b^*$ its mirror
image. Suppose $K_b$ is a knot. We have: 
$$S_b(q,p)=\frac{I(\ra)(K_b^*) I(\ra)(K_b)}{[2\a+1]^2}=X(0,p,K_b)(h)\frac{(2\a+1)^2}{[2a+1]^2},$$
where $q=\exp(h/2)$. The last equality follows from theorem \ref{JonesRelation}.
\end{Proposition}
Notice also lemma \ref{truncate} and that  $X(0,p,K_b)(h)$ is a convergent power series if $p\in\N$. Therefore the perturbative framework of the previous sections is correct, at least for finite dimensional representations. In the sequel we will generalise this for infinite dimensional representations.

\subsubsection{The series Are Convergent $h$-Adicaly}

We now define the $h$-adic version of the theory  developed by Buffenoir and Roche. Let $q \in (0,1)$  and  consider the element $\gai\in \pol$. For any $p\in \C$, we have a balanced representations $\r(p)$ of the Quantum Lorentz Group in $V(p)$. The term 
$$\left < \vb^\ib|\r(p)(\gai)|\vg_{\ig}\right >_q$$ can be seen  a function of
$q$. Due to the fact the building blocks of $\r(p)$ are Clebsch-Gordan coefficients, it express as a sum of square roots of rational functions of $q$, which extend to a  well
defined analytic function in a neighbourhood of $1$. We can see it for example from
(\ref{CGformula}). In addition we have some terms of the form $q^{p\s}, \s \in \Z$, which after putting $q=\exp(h/2)$ define an analytic function of $h$. Therefore 
$$h \mapsto \left < \vb^\ib|\r(p)(\gai)|\vg_{\ig}\right >_{\exp(h/2)}$$ defines a
power series in $h$, uniquely. In particular it follows that if $b$ is a braid then each term of the sum $S_b(\exp(h/2),p)$ defines uniquely a power series in $h$, which converges to the term for $h$ small enough.
\begin{Lemma}
For any $x\in \pol$, the order of: $$h \mapsto \left <
  \vb^\ib|\r(p)(x)|\vg_{\ig}\right >_{\exp(h/2)}, $$ as a power series in $h$,  is bigger or equal to  $|\b-\g|$.
\end{Lemma}
\begin{Proof}
Notice that $\gh^i_j$ sends $\Vg$ to $\Vgm \op \Vg \op \VgM$, in a way such that for $q=1$ the projection $v$ of  $\gh^i_j \vg_{\ig}$ in $\VgM \op\Vgm$ is
zero. We can see this from lemma \ref{lfor}.  In particular $v$ has order bigger or equal to one.  This lemma is thus a trivial consequence of the fact  the elements   $\left \{\gh^i_j,{-1/2\leq i,j \leq 1/2}\right \}$ generate $\pol$ as an algebra.
\end{Proof}

Therefore
\begin{Proposition}
For any braid $b$ whose closure is a knot the infinite sum $S_b(\exp(h/2),p)$ converges in the $h$-adic topology.
\end{Proposition}
\begin{Proof}
Let $b$ be a braid with $n$ crossings and $m+1$ strands. Recall equation (\ref{FormalRT}) and comments after. Due to the way the
$\Xak^i_j$ as well as $G$ act in $V(p)$,  the previous lemma guaranties that the order of $\left <\voo,\prod_{l=1}^{2n+m} T(\aa,\ii,\jj,l)\vo \right >$  as a power series in $h$ is bigger or equal to  $\a_k$, for $k=1,..,n$;
and the result follows.
\end{Proof}
\subsubsection{The Series Define a $\Ch$-Valued Knot Invariant} Since we have  proved the $h$-adic convergence of the sums $S_b(\exp(h/2),p)$ to a formal power series, we could now use Markov's theorem and prove that the assignment  $b \mapsto S_b(\exp(h/2),p)$ defines a knot invariant. However the best way to prove this is to reduce it to the finite dimensional case, since we already know that it defines a knot invariant  and the exact form of it, see proposition \ref{ThePoint}. Consider a coefficient $\La{A}{B}{C}{D}(p)_q$ at $q=\exp(h/2)$, thus it is a power series in $h$ convergent for $h$ small enough. From equation (\ref{DefLambda}), we can see that the
dependence of each term in $p$ is polynomial. In particular:
\begin{Lemma}\label{PolDep}
Let $b$ be a braid, consider the power series $S_b(\exp(h/2),p)$ as a function of $p$, the parameter defining a balanced representation of the Quantum Lorentz
Group. Then each term in the expansion of $S_b(\exp(h/2),p)$ as a power series in $h$ is a polynomial in $p$.
\end{Lemma}
\begin{Proof}
Suppose
$A(p)=\sum_{n\in \N_0} A_n(p)h^n$ and $B(p)=\sum_{n \in \N_0} B_n(p) h^n$ are power series whose coefficients depend polynomially in $p$, for example a power series such as $\exp(m p h/2)$. Then also the coefficients of their product depend polynomially in $p$. This immediately proves this lemma. Note that the Clebsch-Gordan coefficients as well as the actions of $G$ and of the elements $\Xai$ do not depend on $p$.
\end{Proof}

Therefore
\begin{Theorem}\label{equivalence}
Let $p \in \C$ and $b$ be a braid. Let also $K_b$ be the closure of $b$. Let also $K_b^*$ be the mirror image of $K$. We have:
$$S_b(\exp(h/2),p)=\frac{X(0,p,K)(2\a+1)^2}{ [2\a +1]^2},$$
where $\a=(p-1)/2$. 
\end{Theorem}
Recall that by theorem \ref{UNFRAMED} the knot invariant $X(0,p)$ is unframed.

\begin{Proof}
By  lemma \ref{PolDep},  we only need to prove this theorem for $p\in \N$. In this case, if $q=(0,1)$ then $S_b(q,p)$ truncates to a finite sum which from the comments after proposition \ref{ThePoint}  equals $X(0,p,K)(2\a+1)^2 [2\a+1]^{-2}$ at $q=\exp(h/2)$. Recall this power series are convergent if $p$ is integer. Each term of the finite sum  $S_b(\exp(h/2),p)$ is a power series in $h$ convergent for $h$ small enough and coinciding with $S_b(q,p)$, for $q \in (0,1)$ and close enough to $1$; thus the result follows.
\end{Proof}

\section*{Acknowledgements}
This work was realised in the course of my PhD  in the University of Nottingham under the supervision of Dr John W. Barrett. I was financially supported by the programme {\em `` PRAXIS-XXI''}, grant number $SFRH/BD/1004\\ /2000$
of {\em Funda\c{c}\~{a}o para a Ci\^{e}ncia e a Tecnologia} (FCT), financed by the  European Community fund {\em Quadro Comunit\'{a}rio de Apoio III}, and also by {\em Programa Operacional ``Ci\^{e}ncia, Tecnologia, Inova\c{c}\~{a}o''} (POCTI) of the {\em Funda\c{c}\~{a}o para a Ci\^{e}ncia e a Tecnologia } (FCT), cofinanced by the European Community fund FEDER. The last stage of this work was financed by the FCT post-doc grant $SFRH/BDP/17552/2004$, part of the research project $POCIT/MAT/60352/2004$ (``Quantum Topology'').

\end{document}